\documentclass[12pt]{article}
\usepackage[dvipdfm]{graphicx,color}
\definecolor{red}{rgb}{1,0,0}
\usepackage[reqno]{amsmath}
\usepackage{amsthm}
\usepackage{amsfonts}
\usepackage{amssymb}
\usepackage{color} 
\usepackage{amsbsy}

\newtheorem{THEorem}{Theorem}

\newtheorem{cor}{Corollary}[section]

\newtheorem{lem}{Lemma}[section]
\newtheorem{prop}{Proposition}[section]

\newtheorem{thm}{Theorem}[section]

\newfont{\bbf}{msbm10 at 12pt}
\def\CC{{\mbox{\bbf C}}}
\def\DD{{\mbox{\bbf D}}}
\def\PP{{\mbox{\bbf P}}}
\def\NN{{\mbox{\bbf N}}}
\def\QQ{{\mbox{\bbf Q}}}
\def\RR{{\mbox{\bbf R}}}
\def\ZZ{{\mbox{\bbf Z}}}

\def\C{{\cal C}}

\def\qed{$\Box$}
\def\Re{\mbox{\rm Re}}
\def\Im{\mbox{\rm Im}}

\newcommand{\pr}{{\it proof.}\quad}
\newcommand{\ovl}{\overline}

\newcommand{\eps}{\epsilon}

\begin{document}
\centerline{\bf Postcritical sets and saddle basic sets for Axiom A 
polynomial skew products on $\CC^2$}
\vskip 5 mm

\begin{flushright}
{Shizuo Nakane\quad Tokyo Polytechnic University}
\end{flushright}


{\footnotesize  Investigating the link between postcritical behaviors and the 
relations of saddle basic sets for Axiom A polynomial skew products on $\CC^2$, 
we characterize various properties concerning the three kinds of accumulation 
sets defined in \cite{dh1} in terms of the saddle basic sets. We also give 
a new example of higher degree. }
 
\section{Introduction}

In holomorphic dynamics, the behaviors of the orbits of critical points play a 
central role. And so do the relations between the basic sets 
in the dynamics of Axiom A maps. In this article, we investigate the 
link between postcritical behaviors and the relations of saddle basic sets 
for Axiom A polynomial skew products on $\CC^2$. Through this investigation, 
we improve some of the results in \cite{dh1, dh2} and give complete 
characterizations of equalities between three kinds of accumulation sets of 
the critical set. As a corollary, we give stability results of these 
equalities. We also give a new example of higher degree. 

A polynomial skew product on $\CC^2$ is a map of the form : $f(z,w) = (p(z),q(z,w))$, 
where $p(z)$ and $q_z(w) := q(z,w)$ are polynomials of degree $d\geq 2$. 
It is called {\it regular} if it extends to a holomorphic map on $\PP^2$. 
Its $k$-th iterate is expressed by : 
$$
f^k(z,w) = (p^k(z),q_{p^{k-1}(z)}\circ\cdots\circ q_z(w)) =: (p^k(z),Q_z^k(w)). 
$$
Thus it preserves the family of {\it fibers} $\{z\}\times\CC$ and this makes it 
possible to study its dynamics more precisely. Let $K$ be the set of points with 
bounded orbits and set $K_z = \{w\in\CC;(z,w)\in K\}$. 
The {\it fiber Julia set} $J_z$ is the boundary of $K_z$. More generally, 
for $T\subset\CC^2$ and $Z\subset\CC$, we set $T_z = \{w\in\CC;(z,w)\in T\}$ 
and $T_Z = \cup_{z\in Z}\{z\}\times T_z$. 

Let $K_p$ and $J_p$ be the filled-in Julia set and Julia set of $p$ respectively. 
We are interested in the dynamics of $f$ on $J_p\times\CC$ because the dynamics 
outside $J_p\times\CC$ is fairly simple. Consider the {\it critical set} 
$$
C_Z = \{(z,w)\in Z\times\CC; q_z'(w) = 0\}
$$ 
over $Z\subset\CC$. We will investigate the behaviors of the orbits of points 
in $C_{J_p}$. For any subset $X$ in $\CC^2$, its accumulation set is defined by 
$$
A(X) = \cap_{N\geq 0}\ovl{\cup_{n\geq N}f^n(X)}.
$$
DeMarco \& Hruska \cite{dh1} defined the {\it pointwise} and {\it component-wise} 
accumulation sets of $C_{J_p}$ respectively  by
$$
A_{pt}(C_{J_p}) = \ovl{\cup_{x\in C_{J_p}}A(x)}\: \mbox{ and }\:
A_{cc}(C_{J_p}) = \ovl{\cup_{C\in \C(C_{J_p})}A(C)}, 
$$
where $\C(C_{J_p})$ denotes the collection of connected components of $C_{J_p}$. 
It follows from the definition that 
$$
A_{pt}(C_{J_p})\subset A_{cc}(C_{J_p})\subset A(C_{J_p}). 
$$
It also follows that $A_{pt}(C_{J_p}) = A_{cc}(C_{J_p})$ if $J_p$ is a Cantor set 
while $A_{cc}(C_{J_p}) = A(C_{J_p})$ if $J_p$ is connected. 

Let $\Lambda$ be the closure of the set of saddle periodic points in $J_p\times\CC$. 
The {\it stable} and {\it unstable sets} of $\Lambda$ are respectively defined by 
\begin{eqnarray*}
W^s(\Lambda) &=& \{y\in\CC^2;f^k(y)\to\Lambda\},\\
W^u(\Lambda) &=& \{y\in\CC^2;\exists\mbox{ prehistory }\hat{y} 
= (y_k)_{k\leq 0} \mbox{ of $y$ tending to } \Lambda\}.
\end{eqnarray*} 

Throughout this section, we assume that $f$ is an Axiom A regular polynomial 
skew product on $\CC^2$. 
In \cite{dh1}, they characterized $A_{pt}(C_{J_p})$ and $A(C_{J_p})$. 

\begin{THEorem}{\rm (\cite{dh1}, Theorem 1.1)}\label{ThmCh}
$$
A_{pt}(C_{J_p}) = \Lambda, \quad A(C_{J_p}) = W^u(\Lambda)\cap (J_p\times\CC). 
$$
\end{THEorem} 

Thus, any $x\in C_{J_p}$ either tends to $\Lambda$ or escapes to $\infty$. 
They also tried to characterize the equalities between any two of the accumulation 
sets and give examples with various properties. 

We will improve their works and give a new example of higher degree. 
A key tool is the {\it saddle basic set}, 
i.e., the basic set of unstable dimension one. Here 
a {\it basic set} is a compact invariant subset of the non-wandering 
set $\Omega$ of $f$ with dense orbit. 
The saddle set $\Lambda$ decomposes into a disjoint union of 
saddle basic sets : $\Lambda = \sqcup_{i=1}^m\Lambda_i$, 
which is the saddle part of $\Omega$ in $J_p\times\CC$. 
For convenience, we add one more ``basic set,'' which corresponds to the 
superattracting fixed point $\{[0:1:0]\}$ in $\PP^2$ : 
$$
\Lambda_0 = \emptyset,\:\: W^s(\Lambda_0) = (J_p\times\CC)\setminus K,\:\:
W^u(\Lambda_0) = \emptyset. 
$$
Note that this is not a saddle basic set. We also set 
$$
C_i = C_{J_p}\cap W^s(\Lambda_i)\:\: (0\leq i\leq m). 
$$ 

First we characterize the equality $A_{cc}(C_{J_p}) = A_{pt}(C_{J_p})$. 

\begin{thm}\label{ThmAptEqAcc}
$$A_{cc}(C_{J_p}) = A_{pt}(C_{J_p})\Longleftrightarrow \forall C\in\C(C_{J_p}),\: 
0\leq\exists i\leq m\:\mbox{ such that }\: C\subset C_i. \eqno{(1)}$$
\end{thm}

This is an improvement of the following. 

\begin{THEorem}{\rm (\cite{dh1,dh2})}\label{ThmAccAptDH}
$$
A_{cc}(C_{J_p}) = A_{pt}(C_{J_p}) \Longrightarrow \forall C\in\C(C_{J_p}), \:
      C\cap K = \emptyset \mbox{ or } C\subset K. \eqno{(2)}
$$
\end{THEorem}

In fact, in terms of $C_i$, the condition in (2) is expressed by 
$$
\forall C\in\C(C_{J_p}),\quad C\subset C_0 \mbox{ or } C\subset\cup_{i=1}^m C_i.
$$ 
Hence, the condition in (1) coincides with that in (2) only if $m = 1$, 
that is, $\Lambda$ itself is a basic set. 

The following two theorems give characterizations of $A(C_{J_p}) = A_{pt}(C_{J_p})$ 
in terms of $C_i$. 

\begin{thm}\label{ThmAptEqA}
For each $i\geq 0$, 
$$
A(C_i) = \Lambda_i \Longleftrightarrow C_i \mbox{ is closed}.     \eqno{(3)}
$$
Consequently, 
$$
A(C_{J_p}) = A_{pt}(C_{J_p}) \Longleftrightarrow C_i \mbox{ is closed for any } i\geq 0. 
$$
\end{thm}

\begin{thm}\label{ThmLzConti} For each $j\geq 1$,
\begin{eqnarray*}
C_j \mbox{ is open in } C_{J_p} 
&\Longleftrightarrow& W^u(\Lambda_j)\cap (J_p\times\CC) = \Lambda_j\\
&\Longleftrightarrow& z\mapsto \Lambda_{j,z}\mbox{ is continuous in } J_p.
\end{eqnarray*}
Consequently, 
\begin{eqnarray*}
\forall j\geq 1, C_j\mbox{ is open in}\ C_{J_p} 
&\Longleftrightarrow& W^u(\Lambda)\cap (J_p\times\CC) = \Lambda\\ 
&\Longleftrightarrow& z\mapsto\Lambda_z \mbox{ is continuous in}\ J_p. 
\end{eqnarray*}
\end{thm}

Note that $C_0 = C_{J_p}\setminus K$ is always open in $C_{J_p}$. 
Here the continuity is with respect to the Hausdorff topology. 
These can be regarded as a refinement of the following. 

\begin{THEorem}\label{ThmE5.2}{\rm (\cite{dh2}, Theorem E5.2)}
$$
A(C_{J_p}) = A_{pt}(C_{J_p}) \Longleftrightarrow 
      \mbox{ the map } z\mapsto\Lambda_z \mbox{ is continuous in } J_p. \eqno{(4)}
$$
Under the assumption $W^u(\Lambda)\cap W^s(\Lambda) = \Lambda$, 
$$
A(C_{J_p}) = A_{pt}(C_{J_p}) \Longleftrightarrow 
\mbox{ the map } z\mapsto K_z \mbox{ is continuous in } J_p. \eqno{(5)}
$$
\end{THEorem}

Here we give some examples. 

\noindent
{\bf Example 1.1.} 
The product map $f(z,w) = (p(z),q(w))$ is Axiom A if and only if both 
$p$ and $q$ are hyperbolic. Let $A_i = \{w_{i,1},\cdots, w_{i,k_i}\}, 1\leq i\leq m$ 
be the attracting cycles of $q$ and $B_i$ be the set of critical points of $q$ 
contained in the basin of $A_i$. Then $\Lambda_i = J_p\times A_i$ are the saddle basic 
sets and $C_i = J_p\times B_i$ are open and closed in $C_{J_p}$. 

\vspace{2mm}
\noindent
{\bf Example 1.2.} 
Sumi \cite{su} gives an example of a polynomial skew product : 
$$
f(z,w) = \left(p(z), w^{2^n} + \frac{z+\sqrt{R}}{2\sqrt{R}}t_{n,\eps}(w)\right), 
$$
where $R, \eps > 0, n\in\NN, p = p_R^n, p_R(z) = z^2-R, t_{n,\eps}(w) = h_{\eps}^n(w) 
- w^{2^n}, h_{\eps}(w) = (w-\eps)^2-1+\eps.$ If $R$ is large, $\eps$ is small and $n$ 
is even and large, then $f$ is Axiom A and $J_p$ is a Cantor set. Let $\alpha < 0$ 
and $\beta > 0$ be the fixed points of $p_R$. Then it satisfies 

\vspace{2mm}
\noindent
$(a)$ $\Lambda = \Lambda_1\sqcup\Lambda_2$, where $\Lambda_1$ consists of a single point 
in $\{\beta\}\times\CC$, 

\noindent
$(b)$ $C_{J_p}\subset K$ i.e. $C_0 = \emptyset$, hence $z\mapsto K_z$ is continuous 
in $J_p$, 

\noindent
$(c)$ $C_1$ is a finite set in $\{\beta\}\times\CC$, 

\noindent
$(d)$ $C_2 = C_{J_p}\setminus C_1$ is open in $C_{J_p}$ and $\ovl{C_2}\supset C_1$, 

\noindent
$(e)$ $A_{pt}(C_{J_p}) = A_{cc}(C_{J_p}) \neq A(C_{J_p})$. 

\vspace{2mm}
\noindent
By Theorems \ref{ThmAptEqA} and \ref{ThmLzConti}, it follows that $A(C_1) = \Lambda_1$ 
and $W^u(\Lambda_2)\cap (J_p\times\CC) = \Lambda_2$. From $(d)$, we have 
$W^u(\Lambda_1)\cap W^s(\Lambda_2)\neq\emptyset$ (see Proposition \ref{PropACiWs} below).  

\vspace{2mm}
Thus, the equivalence in (5) does not hold in general. So far, Example 1.2 is 
the only example of an Axiom A map which has two saddle basic sets 
with a relation. This suggests that we need to take into account the 
relations among saddle basic sets. 
It turns out that the equality $A(C_{J_p}) = A_{pt}(C_{J_p})$ or 
$W^u(\Lambda)\cap (J_p\times\CC) = \Lambda$ decomposes into two independent equalities : 
$$
W^u(\Lambda)\cap (J_p\times\CC) = W^u(\Lambda)\cap W^s(\Lambda) \mbox{ and } 
W^u(\Lambda)\cap W^s(\Lambda) = \Lambda, 
$$
both of which are characterized in terms of $C_i$. 
See Theorems \ref{ThmKzConti} and \ref{ThmNoRelation}. 
As will be seen in Lemma \ref{LemRelation}, 
$$
W^u(\Lambda)\cap W^s(\Lambda) = \Lambda \Longleftrightarrow 
W^u(\Lambda_i)\cap W^s(\Lambda_j) = \emptyset \mbox{ for any } 1\leq i\neq j\leq m. 
\eqno{(6)}
$$ 
As for the stability of the equalities $A_{cc}(C_{J_p}) = A_{pt}(C_{J_p})$ and 
$A(C_{J_p}) = A_{pt}(C_{J_p})$, we have the following. 

\begin{thm}\label{ThmStable} 
Both equalities $A_{cc}(C_{J_p}) = A_{pt}(C_{J_p})$ and $A(C_{J_p}) = A_{pt}(C_{J_p})$ 
are preserved in hyperbolic components. 
\end{thm} 

The idea of proof of Theorem \ref{ThmStable} is due to that of Proposition 6.3 
in \cite{dh1}. Note that \cite{dh2} has already given another proof for the equality 
$A(C_{J_p}) = A_{pt}(C_{J_p})$, based on the characterization in Theorem \ref{ThmE5.2}. 
By virtue of Theorems \ref{ThmAptEqAcc} and \ref{ThmAptEqA}, we can prove 
both cases in the same way. 

The following two theorems give answers to some questions in \cite{dh1}. 

\begin{thm}\label{ThmAccEqA}
Suppose $J_p$ is disconnected. Then 
$$
A_{cc}(C_{J_p}) = A(C_{J_p})\Longleftrightarrow A_{pt}(C_{J_p}) = A_{cc}(C_{J_p}) 
= A(C_{J_p}). 
$$ 
\end{thm}

Note that $A_{cc}(C_{J_p}) = A(C_{J_p})$ if $J_p$ is connected. Thus this gives a 
characterization of the equality $A_{cc}(C_{J_p}) = A(C_{J_p})$. Together with 
Theorem \ref{ThmStable}, it follows that the equality $A_{cc}(C_{J_p}) = A(C_{J_p})$ 
is preserved in hyperbolic components. This answers Question 8.2 in \cite{dh1}. 

\begin{thm}\label{ThmApt=AccNotA}
There exists an Axiom A polynomial skew product $f$ of degree four with the 
following properties : \\
$(a)$ $J_p$ is neither connected nor totally disconnected, \\
$(b)$ $A_{pt}(C_{J_p}) = A_{cc}(C_{J_p})\neq A(C_{J_p})$. 
\end{thm}

In \cite{dh1}, they give examples satisfying each of the following properties 
respectively except $(vi)$:

\vspace{2mm}
\noindent
$(i)$ $A_{pt}(C_{J_p}) = A_{cc}(C_{J_p}) = A(C_{J_p}),$ 

\noindent
$(ii)$ $A_{pt}(C_{J_p}) \neq A_{cc}(C_{J_p}) = A(C_{J_p}),$ 

\noindent
$(iii)$ $A_{pt}(C_{J_p}) = A_{cc}(C_{J_p}) \neq A(C_{J_p}),$ 

\noindent
$(vi)$ $A_{pt}(C_{J_p}) \neq A_{cc}(C_{J_p}) \neq A(C_{J_p}).$ 

\vspace{2mm}
But all of their base Julia sets $J_p$ are either connected or totally disconnected and 
they posed a question as Question 8.1 on the existence of examples whose base Julia sets 
are neither connected nor totally disconnected. Theorem \ref{ThmAccEqA} says that there 
is no such example satisfying $(ii)$, while Theorem \ref{ThmApt=AccNotA} gives one 
satisfying $(iii)$. It is still unknown whether there exists an example satisfying $(vi)$. 

The author would like to thank Hiroki Sumi for helpful discussion on his example 
and Laura DeMarco for valuable advice. 
He also thanks the referees for helpful suggestions that refine this article.

\section {Preliminaries} 

In this section, first we prepare several notions and properties in the theory of 
dynamics of hyperbolic $C^{\infty}$ endomorphisms. This theory was established 
in Przytycki \cite{p} and Ruelle \cite{r}. Here we collect them from Jonsson \cite{j1}. 
Let $f$ be a $C^{\infty}$ endomorphism of a $C^{\infty}$ Riemannian manifold $M$. 
Consider a compact set $L\subset M$ which satisfies $f(L) = L$. 
The hyperbolicity of $L$ for a non-invertible map $f$ is defined through the 
{\it natural extension} : 
$$ 
\hat{L} = \{\hat{x} = (x_k)_{k\leq 0}; x_j\in L, f(x_k) = x_{k+1}\:\: (k\leq -1)\}. 
$$
and the invertible shift map $\hat{f} : \hat{L}\to\hat{L},\: \hat{f}((x_k)) = (x_{k+1}).$

An endomorphism $f$ is said to be {\it Axiom A} if the {\it non-wandering set} 
$\Omega$ is compact, periodic points are dense in $\Omega$ 
(hence $f(\Omega) = \Omega$) and $\Omega$ is hyperbolic. The following lemma is a 
consequence of the fact that the natural extension of a hyperbolic set for an 
open Axiom A endomorphism has local product structure. 
See Proposition 3.3 in \cite{j1}. 

\begin{lem}\label{LemFundNbd}{\rm (\cite{j1}, Corollary 2.6)}\\
Let $L$ be a hyperbolic set for an open Axiom A endomorphism $f$. 
Then, for any sufficiently small neighborhood $U$ of $L$, we have 

\vspace{2mm}
\noindent
$(a)$ if $y\in U$ and $f^k(y)\in U$ for any $k\geq 0$, then $y\in W_{loc}^s(x)$ 
for some $x\in L$, \\
$(b)$ if $y\in U$ has a prehistory $\hat{y} = (y_k)$ with $y_k\in U$ for 
any $k\leq 0$, then $y\in W_{loc}^u(\hat{x})$ for some $\hat{x}\in\hat{L}$. 
\end{lem}

By Corollary 3.5 in \cite{j1} or Theorem A.3 in \cite{j3}, for an open Axiom A 
endomorphism $f$, the non-wandering set $\Omega$ has a {\it spectral decomposition} 
$\Omega = \sqcup_i\Omega_i$ into basic sets. Here a subset $\Omega_i$ 
of $\Omega$ is called a {\it basic set} if it is compact, 
satisfies $f(\Omega_i) = \Omega_i$ and $f$ is transitive on $\Omega_i$. 
A {\it relation} $\succ$ is defined between basic sets by $\Omega_i\succ\Omega_j$ if 
$(W^u(\Omega_i)\setminus\Omega_i)\cap (W^s(\Omega_j)\setminus\Omega_j)\neq \emptyset$. 
A {\it cycle} is a chain of basic sets satisfying 
$$
\Omega_{i_1}\succ\Omega_{i_2}\succ\cdots\succ\Omega_{i_n} = \Omega_{i_1}. 
$$
There is no trivial cycle for open Axiom A endomorphisms : 

\begin{lem}\label{LemNoTrivialCycle}
{\rm (\cite{j1}, Lemma 4.1 or \cite{j3}, Proposition A.4)}\\ 
For open Axiom A endomorphisms, 
$W^u(\Omega_i)\cap W^s(\Omega_i) = \Omega_i$ holds for any $i$. 
\end{lem}

In case $i\neq j$, the property $\Omega_i\succ\Omega_j$ is somewhat simplified. 
The following also shows the equivalence in (6). 

\begin{lem}\label{LemRelation}
For $i\neq j$, $\Omega_i\succ\Omega_j$ if and only if 
$W^u(\Omega_i)\cap W^s(\Omega_j)\neq\emptyset$. 
\end{lem}
\pr  We have only to show the `if' part. 
Suppose $W^u(\Omega_i)\cap W^s(\Omega_j)\neq\emptyset$ but 
$(W^u(\Omega_i)\setminus\Omega_i)\cap (W^s(\Omega_j)\setminus\Omega_j) = \emptyset$. 
Since $W^u(\Omega_i)\cap W^s(\Omega_j)$ never intersects $\Omega_i$, 
it is included in $\Omega_j$. Any point $x\in W^u(\Omega_i)\cap W^s(\Omega_j)$ has 
a prehistory $\hat{x} = (x_k)$ tending to $\Omega_i$. Then, for all $k\leq 0$, 
we have $x_k\in W^u(\Omega_i)\cap W^s(\Omega_j)$, hence $x_k\in\Omega_j$. 
This is a contradiction. 
\qed

\vspace{2mm}
Now we restrict our maps to regular polynomial skew products 
$f(z,w) = (p(z),q(z,w))$ on $\CC^2$ and give some of their basic properties 
from \cite{j3}, which will be repeatedly used later. 
See also Sester \cite{se} for hyperbolicity of fibered polynomials. 

Let $Z$ be a closed subset of $\CC$ such that $p(Z)\subset Z$. Actually, 
$Z$ is either $J_p$ or the set $A_p$ of attracting periodic points of $p$. 
Let $D_Z = \ovl{\cup_{n\geq 1}f^n(C_Z)}$ be the {\it postcritical set} of $C_Z$ 
and set $J_Z = \ovl{\cup_{z\in Z} \{z\}\times J_z}$. 
Jonsson \cite{j3} gave a characterization for $f$ to be Axiom A. 

\begin{thm}\label{ThmAxiomA}{\rm (\cite{j3}, Theorems 8.2 and 3.1)}
A regular polynomial skew product $f$ on $\CC^2$ is Axiom A if and only if 
the following three conditions are satisfied : 

\vspace{2mm}
\noindent
$(a)$ $p$ is hyperbolic, \\
$(b)$ $D_{J_p}\cap J_{J_p} = \emptyset$, \\
$(c)$ $D_{A_p}\cap J_{A_p} = \emptyset$.
\end{thm}

\begin{thm}\label{ThmVertExp}{\rm (\cite{j3}, Proposition 3.5)}\\
If $D_Z\cap J_Z = \emptyset$, then the map $z\mapsto J_z$ is 
continuous in $Z$. 
\end{thm}

Let $\mu$ be the ergodic measure of maximal entropy for $f$ (see \cite{fs1} or \cite{j3}). 
Its support $J_2$ is called the {\it second Julia set} of $f$. 
Corollary 4.4 in \cite{j3} says that $J_2 = J_{J_p}$. By Theorems \ref{ThmAxiomA} 
and \ref{ThmVertExp}, $J_2 = \cup_{z\in J_p}\{z\}\times J_z$ if $f$ is Axiom A. 

The following is a key lemma for the proof of Theorem \ref{ThmAptEqAcc}. 

\begin{thm}{\rm (\cite{j3}, Theorem 8.2 and Proposition A.7)}\label{ThmNoCycle}\\
Axiom A regular polynomial skew products on $\CC^2$ have no cycles.
\end{thm}

Note that the local stable manifold $W_{loc}^s(x)$ of $x\in\Lambda$ is included 
in the fiber containing $x$ because $f$ is contracting in the fiber direction 
on $\Lambda$. Hence the local unstable manifold $W_{loc}^u(\hat{x})$ 
is transversal to the fiber for any $\hat{x}\in\hat\Lambda$. 

We remark that Fornaess \& Sibony \cite{fs2} also investigated 
hyperbolic holomorphic maps on $\PP^2$. 
See also Mihailescu \cite{mih} and Mihailescu \& Urbanski \cite{mu}.

\section{Proofs of Theorems}
\subsection{Proof of Theorem \ref{ThmAptEqAcc}}

Note that $A_{cc}(C_{J_p}) = A_{pt}(C_{J_p})$ if and only if $A(C)\subset\Lambda$ 
for any $C\in\C(C_{J_p})$. First we show the following. 

\begin{lem}\label{LemCritSet} 
For Axiom A regular polynomial skew product $f$ on $\CC^2$, we have 
$C_{J_p} = \sqcup_{i=0}^m C_i$.
\end{lem}
\pr  Theorem \ref{ThmCh} implies $A(x)\subset\Lambda$ for any $x\in C_{J_p}$. 
If $A(x) = \emptyset$, then $x\in C_0$. Otherwise, by Lemma \ref{LemFundNbd}, 
there exist an $n$ and $y\in\Lambda$ such that $f^n(x)\in W^s_{loc}(y)$. 
Hence $A(x)\subset\Lambda_i$ if $y\in\Lambda_i$. 
Thus we have $C_{J_p} = \sqcup_{i=0}^m C_i$. 
\qed

\vspace{2mm}
Now we prove Theorem \ref{ThmAptEqAcc}. 

\vspace{2mm}
\noindent
($\Rightarrow$)\:  Suppose $C\in\C(C_{J_p})$ intersects at least two of $C_i$. 
By Theorem \ref{ThmAccAptDH}, we may assume $C\subset \cup_{i=1}^m C_i$. 
Then, by Lemma \ref{LemCritSet}, we have $C = \sqcup_{i=1}^m (C\cap C_i)$. 
Since $C$ is closed, if all $C\cap C_i$ are closed, it contradicts the 
connectivity of $C$. Thus at least one of them is not closed. 
We may assume $C\cap \ovl{C_i}\cap C_j\neq\emptyset$ for some $i\neq j$. 
The following holds for $i, j\geq 0$ and will also be used later. 

\begin{lem}\label{LemAccum}
Suppose $i\neq j$ and $\ovl{C_i}\cap C_j\neq\emptyset$. Then 
$A(C_i)\cap (W^u(\Lambda_j)\setminus\Lambda)\neq\emptyset$. 
\end{lem}
\pr  Since $\ovl{C_i}\cap C_0 = \emptyset$ for $i\neq 0$, we may assume $j\geq 1$. 
Take a sequence $x_n\in C_i$ tending to $x_0\in C_j$. 
Take a small open neighborhood $U_k$ of $\Lambda_k$ for $1\leq k\leq m$ 
so that $f(U_k)\cap U_{\ell} = U_k\cap U_{\ell} = \emptyset$ for $k\neq\ell$. 
Since $x_0\in C_j$, there exists a $K\geq 0$ such that $f^k(x_0)\in U_j$ for $k\geq K$. 
Then $f^K(x_n)\in U_j$ for large $n$. Since $x_n\in C_i$, the orbit of $x_n$ 
eventually leaves $U_j$. Hence if we set $k_n = min\{k\geq K; f^k(x_n)\notin U_j\}$, 
then $k_n\to\infty$. 
Let $y$ be an accumulation point of the sequence $\{f^{k_n}(x_n)\}$. 
Then $y\in\ovl{f(U_j)}\setminus U_j$ since $f^{k_n-1}(x_n)\in U_j$. 
Consequently we have $y\notin\cup U_k$, hence $y\in A(C_i)\setminus\Lambda$. 

Next we show $y\in W^u(\Lambda_j)$. Taking subsequences if necessary, 
set $y_{\ell} = \lim_{n\to\infty} f^{k_n+\ell}(x_n)$ for $\ell\leq 0$. 
Then $(y_{\ell})_{\ell\leq 0}$ is a prehistory of $y$ with $y_{\ell}\in\ovl{U_j}$ 
for $\ell\leq -1$. By Lemma \ref{LemFundNbd}, $y_{-1}\in W_{loc}^u(\hat{x})$ 
for some $\hat{x}\in\hat\Lambda_j$, hence $y\in W^u(\Lambda_j)$. 
\qed 

\vspace{2mm}
In the above proof, if we take $x_n\in C\cap C_i$, then $x_0\in C\cap C_j$ 
and $y\in A(C)$. Thus $A(C)$ contains a point $y$ 
outside $\Lambda = A_{pt}(C_{J_p})$. 
Now we conclude $A_{cc}(C_{J_p})\neq A_{pt}(C_{J_p})$. 

\vspace{2mm}
\noindent
($\Leftarrow$)\:  We have only to show that $A(C)\subset\Lambda_i$ if $C\in\C(C_{J_p})$ 
satisfies $C\subset C_i$. More generally, we show the following. 

\begin{lem}\label{LemClosedAccum}
If $C\subset C_i$ is closed, then $A(C)\subset\Lambda_i$. 
\end{lem}
\pr  If $C\subset C_0$, then $A(C) = \emptyset$ since $C$ is compact. 
Suppose $C\subset C_i$ and there exists $x\in A(C)\setminus\Lambda_i$ for $i\geq 1$. 
By taking $U_i$ small, there exist $m_n\nearrow\infty$ and $x_n\in C$ 
such that $f^{m_n}(x_n)\notin U_i$ for any $n$. 
Since $C$ is closed, we may assume $x_n$ tends to some $x_0\in C\subset C_i$. 
Set $k_n = min\{k\geq K; f^k(x_n)\notin U_i\}$ as above and
take an accumulation point $y$ of $\{f^{k_n}(x_n)\}$. 
By the above argument, $y\in W^u(\Lambda_i)\setminus\Lambda_i$, 
hence $y\notin W^s(\Lambda_i)$ because of Lemma \ref{LemNoTrivialCycle}. 
Since $y\in A(C)$, it follows that $y\in K_{J_p}\setminus J_2,$ 
which is, by Lemma 3.6 in \cite{dh1}, equal to $W^s(\Lambda)$. 
Thus $y$ must belong to $W^s(\Lambda_{i_1})$ for some $i_1\neq i$. 
That is, we have a sequence $\{f^{k_n}(x_n)\}$ 
in $W^s(\Lambda_i)$ tending to $y\in W^u(\Lambda_i)\cap W^s(\Lambda_{i_1})$, 
hence $\Lambda_i\succ\Lambda_{i_1}$ by Lemma \ref{LemRelation}. 

By successively applying this argument, we have a chain of saddle basic sets :
$$
\Lambda_i\succ\Lambda_{i_1}\succ\Lambda_{i_2}\succ\cdots,\quad i\neq i_1\neq i_2\neq\cdots.
$$
Since there exist only finitely many basic sets, we must have a cycle of them, 
which contradicts Theorem \ref{ThmNoCycle}. 
This completes the proof. 
\qed

\subsection{Proofs of Theorems \ref{ThmAptEqA} and \ref{ThmLzConti}}

Note that Theorem \ref{ThmAptEqA} follows from Lemmas \ref{LemAccum} and 
\ref{LemClosedAccum}. Here we use the following proposition, 
which completely characterizes the existence of relations between saddle 
basic sets in terms of $C_i$. 
Both Theorems \ref{ThmAptEqA} and \ref{ThmLzConti} easily follow 
from it. Set $I = \{0,1,2,\cdots,m\}$. 

\begin{prop}\label{PropACiWs}
For any $i,j\in I$ with $i\neq j$, the following four conditions are equivalent 
to each other. 

\vspace{2mm}
\noindent
$(a)$ $\ovl{C_i}\cap C_j\neq\emptyset,$\\ 
$(b)$ $A(C_i)\cap W^s(\Lambda_j)\neq\emptyset,$\\
$(c)$ $\ovl{W^s(\Lambda_i)}\cap W^s(\Lambda_j)\neq\emptyset,$ \\
$(d)$ $W^s(\Lambda_i)\cap W^u(\Lambda_j)\neq\emptyset.$
\end{prop}

This proposition holds also for $i = j\geq 1$. 
It holds also for $i = j = 0$, if we set $\Lambda_0 = W^u(\Lambda_0) = \{[0:1:0]\}$. 

\pr  Since, for any $i\geq 1$, all the sets 
$\ovl{C_i}\cap C_0, A(C_i)\cap W^s(\Lambda_0), \ovl{W^s(\Lambda_i)}\cap W^s(\Lambda_0)$ 
and $W^s(\Lambda_i)\cap W^u(\Lambda_0)$ are empty, we may assume $j\geq 1$. 

\vspace{2mm}
\noindent
$(a)\Rightarrow (b)$\:  Suppose $\ovl{C_i}\cap C_j\neq\emptyset$. 
Take a sequence $x_n\in C_i$ tending to $x_0\in C_j$. By Lemma \ref{LemFundNbd}, 
there exist $y\in\Lambda_j$ and $k$ such that $f^k(x_0)\in W^s_{loc}(y)$. 
Hence, for any $n > 0$, there exists $L_n$ such that 
$d(f^{\ell}(x_0),f^{\ell-k}(y)) < 1/n$ for $\ell\geq L_n$. For each fixed $n$, 
take $k_n$ so that $\displaystyle d(f^{L_n}(x_{k_n}), f^{L_n}(x_0)) < 1/n$. 
Then it follows $\displaystyle d(f^{L_n}(x_{k_n}),f^{L_n-k}(y)) < 2/n$. 
Since $f^{L_n-k}(y)\in\Lambda_j$, we conclude that $A(C_i)\cap\Lambda_j\neq\emptyset$. 

\vspace{2mm}
\noindent
$(b)\Rightarrow (c)$\:  It is evident since $A(C_i)\subset\ovl{W^s(\Lambda_i)}$. 

\begin{figure}[h]
\centering{
\unitlength 0.1in
\begin{picture}( 34.6500, 16.3000)(  0.9000,-22.3000)
%
\special{pn 8}%
\special{sh 1}%
\special{ar 890 1166 10 10 0  6.28318530717959E+0000}%
\special{sh 1}%
\special{ar 890 1166 10 10 0  6.28318530717959E+0000}%
\special{sh 1}%
\special{ar 890 1166 10 10 0  6.28318530717959E+0000}%
\special{sh 1}%
\special{ar 890 1166 10 10 0  6.28318530717959E+0000}%
\special{sh 1}%
\special{ar 890 1166 10 10 0  6.28318530717959E+0000}%
\special{sh 1}%
\special{ar 878 1166 10 10 0  6.28318530717959E+0000}%
\special{sh 1}%
\special{ar 462 890 10 10 0  6.28318530717959E+0000}%
\special{sh 1}%
\special{ar 462 890 10 10 0  6.28318530717959E+0000}%
\put(0.9000,-23.8500){\makebox(0,0)[lb]{$W^s(\Lambda_i)$}}%
\put(7.8500,-24.0000){\makebox(0,0)[lb]{$W^s(\Lambda_j)$}}%
\put(9.3600,-12.3900){\makebox(0,0)[lb]{$x_0$}}%
\put(2.4100,-9.3400){\makebox(0,0)[lb]{$x_n$}}%
%
\special{pn 8}%
\special{pa 496 920}%
\special{pa 866 1152}%
\special{dt 0.045}%
\special{sh 1}%
\special{pa 866 1152}%
\special{pa 820 1100}%
\special{pa 822 1124}%
\special{pa 800 1134}%
\special{pa 866 1152}%
\special{fp}%
%
\special{pn 8}%
\special{sh 1}%
\special{ar 3346 1456 10 10 0  6.28318530717959E+0000}%
\special{sh 1}%
\special{ar 3346 1442 10 10 0  6.28318530717959E+0000}%
\special{sh 1}%
\special{ar 3346 1442 10 10 0  6.28318530717959E+0000}%
\special{sh 1}%
\special{ar 3336 1442 10 10 0  6.28318530717959E+0000}%
\special{sh 1}%
\special{ar 3336 1428 10 10 0  6.28318530717959E+0000}%
\put(0.9000,-23.8500){\makebox(0,0)[lb]{$W^s(\Lambda_i)$}}%
\put(24.5400,-24.0000){\makebox(0,0)[lb]{$W^s(\Lambda_i)$}}%
\put(31.8400,-24.0000){\makebox(0,0)[lb]{$W^s(\Lambda_j)$}}%
\put(34.1600,-14.5600){\makebox(0,0)[lb]{$y\in\Lambda_j$}}%
\put(21.6400,-10.7900){\makebox(0,0)[lb]{$f^{k+L}(x_n)$}}%
\put(15.1500,-12.8200){\makebox(0,0)[lb]{$f^{k+L}$}}%
\put(29.4100,-21.1000){\makebox(0,0)[lb]{$C_{J_p}$}}%
%
\special{pn 8}%
\special{sh 1}%
\special{ar 3358 1776 10 10 0  6.28318530717959E+0000}%
\special{sh 1}%
\special{ar 3346 1776 10 10 0  6.28318530717959E+0000}%
%
\special{pn 8}%
\special{pa 1330 1340}%
\special{pa 2002 1442}%
\special{fp}%
\special{sh 1}%
\special{pa 2002 1442}%
\special{pa 1940 1412}%
\special{pa 1950 1434}%
\special{pa 1934 1452}%
\special{pa 2002 1442}%
\special{fp}%
\put(26.3900,-19.0600){\makebox(0,0)[lb]{$C_i$}}%
\put(34.1600,-19.3500){\makebox(0,0)[lb]{$C_j$}}%
%
\special{pn 8}%
\special{pa 2884 1036}%
\special{pa 3324 948}%
\special{dt 0.045}%
\special{sh 1}%
\special{pa 3324 948}%
\special{pa 3254 942}%
\special{pa 3272 958}%
\special{pa 3262 982}%
\special{pa 3324 948}%
\special{fp}%
%
\special{pn 8}%
\special{sh 1}%
\special{ar 3346 934 10 10 0  6.28318530717959E+0000}%
\special{sh 1}%
\special{ar 3336 920 10 10 0  6.28318530717959E+0000}%
\put(33.9300,-9.7700){\makebox(0,0)[lb]{$f^{k+L}(x_0)$}}%
%
\special{pn 8}%
\special{sh 1}%
\special{ar 2884 1806 10 10 0  6.28318530717959E+0000}%
\special{sh 1}%
\special{ar 2884 1806 10 10 0  6.28318530717959E+0000}%
%
\special{pn 8}%
\special{pa 450 600}%
\special{pa 450 2082}%
\special{fp}%
%
\special{pn 8}%
\special{pa 890 600}%
\special{pa 890 2066}%
\special{fp}%
%
\special{pn 8}%
\special{pa 3346 616}%
\special{pa 3346 2066}%
\special{fp}%
%
\special{pn 8}%
\special{pa 2906 1806}%
\special{pa 3312 1790}%
\special{dt 0.045}%
\special{sh 1}%
\special{pa 3312 1790}%
\special{pa 3246 1772}%
\special{pa 3260 1792}%
\special{pa 3246 1812}%
\special{pa 3312 1790}%
\special{fp}%
%
\special{pn 8}%
\special{pa 2780 1442}%
\special{pa 2812 1434}%
\special{pa 2842 1426}%
\special{pa 2874 1416}%
\special{pa 2906 1410}%
\special{pa 2938 1402}%
\special{pa 2968 1396}%
\special{pa 3000 1392}%
\special{pa 3032 1388}%
\special{pa 3062 1386}%
\special{pa 3094 1384}%
\special{pa 3124 1386}%
\special{pa 3156 1388}%
\special{pa 3186 1392}%
\special{pa 3218 1398}%
\special{pa 3248 1404}%
\special{pa 3280 1412}%
\special{pa 3310 1422}%
\special{pa 3342 1432}%
\special{pa 3372 1442}%
\special{pa 3402 1454}%
\special{pa 3434 1466}%
\special{pa 3464 1478}%
\special{pa 3494 1490}%
\special{pa 3526 1502}%
\special{pa 3556 1516}%
\special{sp 0.070}%
\put(21.7600,-15.0000){\makebox(0,0)[lb]{$W^u_{loc}(\hat{y})$}}%
%
\special{pn 8}%
\special{pa 3346 1166}%
\special{pa 3346 1298}%
\special{fp}%
\special{sh 1}%
\special{pa 3346 1298}%
\special{pa 3366 1230}%
\special{pa 3346 1244}%
\special{pa 3326 1230}%
\special{pa 3346 1298}%
\special{fp}%
%
\special{pn 8}%
\special{pa 3346 1690}%
\special{pa 3346 1588}%
\special{fp}%
\special{sh 1}%
\special{pa 3346 1588}%
\special{pa 3326 1654}%
\special{pa 3346 1640}%
\special{pa 3366 1654}%
\special{pa 3346 1588}%
\special{fp}%
%
\special{pn 8}%
\special{pa 3208 1398}%
\special{pa 3150 1384}%
\special{fp}%
\special{sh 1}%
\special{pa 3150 1384}%
\special{pa 3210 1420}%
\special{pa 3202 1398}%
\special{pa 3218 1380}%
\special{pa 3150 1384}%
\special{fp}%
%
\special{pn 8}%
\special{pa 3452 1486}%
\special{pa 3544 1516}%
\special{fp}%
\special{sh 1}%
\special{pa 3544 1516}%
\special{pa 3486 1476}%
\special{pa 3492 1498}%
\special{pa 3474 1514}%
\special{pa 3544 1516}%
\special{fp}%
%
\special{pn 8}%
\special{pa 2884 630}%
\special{pa 2884 2082}%
\special{fp}%
%
\special{pn 8}%
\special{sh 1}%
\special{ar 2884 1050 10 10 0  6.28318530717959E+0000}%
\special{sh 1}%
\special{ar 2884 1050 10 10 0  6.28318530717959E+0000}%
\special{sh 1}%
\special{ar 2884 1050 10 10 0  6.28318530717959E+0000}%
\special{sh 1}%
\special{ar 2884 1050 10 10 0  6.28318530717959E+0000}%
\end{picture}%
}
\caption{Local stable and unstable manifolds}
\label{FigStableMfd}
\end{figure}
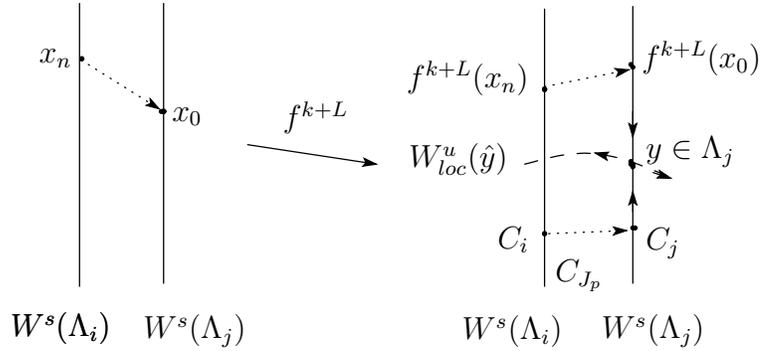

\vspace{2mm}
\noindent
$(c)\Rightarrow (a)$\:  Suppose $W^s(\Lambda_i)\ni x_n\to x_0\in W^s(\Lambda_j)$. 
Then there exist $k$ and $q\in\Lambda_j$ such that 
$f^k(x_0)\in W^s_{loc}(q)$. Let $U_q$ be the connected component of the vertical 
slice of $(J_p\times\CC)\setminus J_2$ containing $q$. 
Then $U_q\supset W^s_{loc}(q)$ and, by Proposition 3.8 in \cite{dh1}, 
there exists $L > 0$ such that $f^L(U_q) = U_{f^L(q)}$ contains a critical point $c$, 
which is in $C_j$. Set $y = f^L(q)$. By Theorems \ref{ThmAxiomA} and \ref{ThmVertExp}, 
any compact subset in $U_y = f^{k+L}(U_{x_0})$ is approximated by compact subsets 
in $f^{k+L}(U_{x_n})$. Hence the branch of the critical locus 
through $c$ must intersect $f^{k+L}(U_{x_n})$ for large $n$. 
Thus, for large $n$, there exist critical points $c_n\in f^{k+L}(U_{x_n})$ tending 
to $c$. Since $c_n\in C_i$, we conclude that $\ovl{C_i}\cap C_j\neq\emptyset.$ 
See Figure \ref{FigStableMfd}. 

\vspace{2mm}
\noindent
$(c)\Rightarrow (d)$\:  Suppose $W^s(\Lambda_i)\ni x_n\to x_0\in W^s(\Lambda_j)$. 
By the above argument, any compact subset in $U_y = U_{f^{k+L}(x_0)}$ is 
approximated by compact subsets in $U_{f^{k+L}(x_n)}$. Since, 
for any prehistory $\hat{y}\in\hat\Lambda_j$ of $y$, $W^u_{loc}(\hat{y})$ is 
transversal to the fiber, it follows that 
$U_{f^{k+L}(x_n)}\cap W^u_{loc}(\hat{y})\neq\emptyset$ for large $n$. 
Thus we conclude $W^s(\Lambda_i)\cap W^u(\Lambda_j)\neq\emptyset.$ 

\vspace{2mm}
\noindent
$(d)\Rightarrow (c)$\:  We need a lemma. 

\begin{lem}\label{LemLocUnstableMfd}
Suppose $\ovl{W^s(\Lambda_i)}\cap W^s(\Lambda_j) = \emptyset.$ Then there exists 
$\delta > 0$ such that $W^s(\Lambda_i)\cap W^u_{\delta}(\hat{y}) = \emptyset$ 
for any $\hat{y}\in\hat\Lambda_j$. 
\end{lem}
\pr  Suppose, for any $n\geq 1$, there exist a prehistory $\hat{y}_n\in\hat\Lambda_j$ 
of $y_n\in\Lambda_j$ and $x_n\in W^s(\Lambda_i)\cap W^u_{1/n}(\hat{y}_n)$. 
Since the sequences $\{x_n\}$ and $\{y_n\}$ are bounded,
there exist their convergent subsequences $\{x_{n_k}\}$ and $\{y_{n_k}\}$ 
respectively tending to some points $x_0$ and $y_0.$ Since $d(x_n,y_n) < 1/n$, 
we have $x_0 = y_0\in\Lambda_j$. 
Thus we conclude $\ovl{W^s(\Lambda_i)}\cap W^s(\Lambda_j)\neq\emptyset.$ 
\qed

\vspace{2mm}
Now suppose $W^s(\Lambda_i)\cap W^u(\Lambda_j)\neq\emptyset.$ 
If we take $p\in W^s(\Lambda_i)\cap W^u(\Lambda_j)$, 
there exists a prehistory $\hat{p} = (p_k)$ of $p$ tending to $\Lambda_j$. 
Assume $\ovl{W^s(\Lambda_i)}\cap W^s(\Lambda_j) = \emptyset.$ Let $\delta > 0$ 
be the constant in Lemma \ref{LemLocUnstableMfd}. Then there exists 
$L > 0$ such that $d(p_k,\Lambda_j) < \delta$ for any $k\leq -L$, 
hence by Lemma \ref{LemFundNbd}, $p_{-L}\in W^u_{\delta}(\hat{y})$ for 
some $\hat{y}\in\hat\Lambda_j$. Since $p_{-L}\in W^s(\Lambda_i)$, 
we have $W^s(\Lambda_i)\cap W^u_{\delta}(\hat{y})\neq\emptyset,$ 
which contradicts Lemma \ref{LemLocUnstableMfd}. 
This completes the proof of Proposition \ref{PropACiWs}. \qed

\vspace{2mm}
Set $I_i = \{j\in I; \ovl{C_i}\cap C_j\neq\emptyset\}$. 
\begin{cor}\label{CorIncStable}
$
A(C_i) \subset \cup_{j\in I_i} W^s(\Lambda_j)\quad  (0\leq i\leq m). 
$
\end{cor}
\pr  By Theorem \ref{ThmCh}, for $i\geq 1$ 
it follows $A(C_i)\subset W^u(\Lambda)\cap K_{J_p}$, which is, 
by Lemma E3.5 in \cite{dh2}, included in $W^s(\Lambda)$. In case $i = 0$, 
$A(C_0) = (A(C_0)\cap W^s(\Lambda_0))\sqcup (A(C_0)\cap W^s(\Lambda))$ 
and $A(C_0)\cap W^s(\Lambda)$ is treated similarly. 
Now the assertion follows from Proposition \ref{PropACiWs}. \qed

\vspace{2mm}
{\it proof of Theorem \ref{ThmAptEqA}.}\:  
Note that $C_i$ is closed if and only if $\ovl{C_i}\cap C_j = \emptyset$ 
for any $j\neq i$. By Proposition \ref{PropACiWs}, this is equivalent to 
$W^s(\Lambda_i)\cap W^u(\Lambda) = W^s(\Lambda_i)\cap W^u(\Lambda_i)$, 
which is, by Lemma \ref{LemNoTrivialCycle}, equal to $\Lambda_i$. 
By Theorem \ref{ThmCh} and Proposition \ref{PropACiWs}, it is also equivalent 
to $A(C_i)\subset W^s(\Lambda_i)\cap W^u(\Lambda) = \Lambda_i$. 
The inclusion $A(C_i)\supset\Lambda_i$ is trivial. This proves (3). 

If $C_i$ for some $i$ is not closed, there exists $j\neq i$ such that 
$\ovl{C_i}\cap C_j\neq \emptyset$. 
By Lemma \ref{LemAccum}, we have $A(C_i)\setminus\Lambda\neq\emptyset$, 
hence $A(C_{J_p})\neq A_{pt}(C_{J_p})$. If $C_i$ is closed for any $i$, then 
$A(C_{J_p}) = A_{pt}(C_{J_p})$ follows from (3). 
\qed

\vspace{2mm}
{\it proof of Theorem \ref{ThmLzConti}.}\:  $C_j$ is open in $C_{J_p}$ if and 
only if $\ovl{C_i}\cap C_j = \emptyset$ for any $i\neq j$. By Proposition 
\ref{PropACiWs} and Lemma \ref{LemNoTrivialCycle}, it is equivalent to 
$W^u(\Lambda_j)\cap (J_p\times\CC) = W^u(\Lambda_j)\cap W^s(\Lambda_j) = \Lambda_j$. 
The equivalence of $W^u(\Lambda_j)\cap (J_p\times\CC) = \Lambda_j$ and 
the continuity of the map $z\mapsto\Lambda_{j,z}$ in $J_p$ can be proved 
by the same way as Theorem \ref{ThmE5.2}. Continuity of the map $z\mapsto\Lambda_z$ 
is equivalent to that of the map $z\mapsto\Lambda_{j,z}$ for any $j\geq 1$. 
\qed

\subsection{Continuity of the map $z\mapsto K_z$}

By virtue of Proposition \ref{PropACiWs}, we will characterize the properties in (5) 
and (6). 

\begin{thm}\label{ThmKzConti} The following three properties are equivalent to each other. \\
$(a)$ $C_0$ is closed, \\
$(b)$ the map $z\mapsto K_z$ is continuous in $J_p$,\\
$(c)$ $A(C_{J_p}) = W^u(\Lambda)\cap W^s(\Lambda)$.
\end{thm}

Recall that the continuity of a set-valued map decomposes into lower and upper 
semicontinuities. 
Theorem \ref{ThmKzConti} reproves the equivalence (5) in Theorem C. 

\pr  Note that $C_0$ is closed if and only if 
$\ovl{C_0}\cap (\cup_{j=1}^m C_j) = \emptyset$, which is, by Proposition 
\ref{PropACiWs}, equivalent to $\ovl{W^s(\Lambda_0)}\cap W^s(\Lambda) = \emptyset$ 
or to $W^s(\Lambda_0)\cap W^u(\Lambda) = \emptyset$. 

\vspace{2mm}
\noindent
$(a)\Leftrightarrow (b)$\: 
Proposition 2.1 in \cite{j3} says that $z\mapsto K_z$ is upper-semicontinuous 
and by Theorem \ref{ThmVertExp}, $z\mapsto J_z$ is continuous in $J_p$. 
Hence, $z\mapsto K_z$ is discontinuous at $z_0\in J_p$ 
if and only if there exist a sequence $z_n$ in $J_p$ tending to $z_0$ 
and $w_0\in int\, K_{z_0}$ such that $w_0\notin K_{z_n}$ for any $n\geq 1$. 
That is,  by Lemma 3.6 in \cite{dh1}, there exists a sequence $(z_n,w_0)$ 
in $W^s(\Lambda_0)$ tending to a point $(z_0,w_0)\in W^s(\Lambda)$. 
This is equivalent to $\ovl{W^s(\Lambda_0)}\cap W^s(\Lambda)\neq\emptyset$, 
that is, $C_0$ is not closed. 

\vspace{2mm}
\noindent
$(a)\Leftrightarrow (c)$\: 
By Theorem \ref{ThmCh} and Lemma E3.5 in \cite{dh2}, we have 
$$
A(C_{J_p}) = (W^u(\Lambda)\cap W^s(\Lambda))\sqcup (W^u(\Lambda)\cap W^s(\Lambda_0)). 
$$ 
Thus $C_0$ is closed if and only if $A(C_{J_p}) = W^u(\Lambda)\cap W^s(\Lambda)$. 
\qed

\vspace{2mm}
The following is a direct consequence of Proposition \ref{PropACiWs} and (6). 
\begin{thm}\label{ThmNoRelation}
$
W^u(\Lambda)\cap W^s(\Lambda) = \Lambda \Longleftrightarrow C_i 
\mbox{ is closed for any } i\geq 1.
$
\end{thm}

\subsection{The sets $A(C_i)$ in general case}

We have shown $A(C_i) = \Lambda_i$ if $C_i$ is closed. In this subsection, 
we extend it and give a description of $A(C_i)$ in general case. 

\begin{thm}\label{ThmEstimate}
Recall that $I_i = \{j\in I; \ovl{C_i}\cap C_j\neq\emptyset\}$. Then we have 
$$
A(C_i)\subset (\cup_{j\in I_i}W^u(\Lambda_j))\cap (\cup_{j\in I_i}W^s(\Lambda_j)).
$$
\end{thm}

If $C_i$ is closed, $I_i = \{i\}$. Then Theorem \ref{ThmEstimate} says 
$A(C_i)\subset W^u(\Lambda_i)\cap W^s(\Lambda_i) = \Lambda_i$. 
Thus Theorem \ref{ThmEstimate} generalizes Theorem \ref{ThmAptEqA}. 

\pr  By virtue of Corollary \ref{CorIncStable}, we have only to show \\
$A(C_i)\subset\cup_{j\in I_i} \left(W^u(\Lambda_j)\cap (J_p\times\CC)\right)$. 
We give a proof mainly for the case $i\geq 1$. 
The case $i = 0$ can be done with a minor change and will be given at the end. 
Below, we repeatedly use the argument in the proof of Lemmas \ref{LemAccum} 
and \ref{LemClosedAccum}. 

Suppose $p\in A(C_i)$. Then there exists $x_n\in C_i$ 
and $m_n\nearrow\infty$ such that $f^{m_n}(x_n)\to p$. 
If $p\in\Lambda$, by Corollary \ref{CorIncStable}, we have $j\in I_i$ 
and $p\in\cup_{j\in I_i}\Lambda_j$. This holds also for $i = 0$. 

In the sequel, we assume $p\notin \Lambda$. 
We may assume $p\notin\ovl{\cup_{j=1}^m U_j}$, where $U_j$ is a small open 
neighborhood of $\Lambda_j$. We may also assume $x_n\to x_0$. If $x_0\in C_i$, 
by Lemma \ref{LemClosedAccum}, we have $p\in\Lambda_i$, which contradicts 
$p\notin \Lambda$. Hence $x_0\in C_{i_1}$ for some $i_1\neq i$ and 
by Proposition \ref{PropACiWs}, $i_1\in I_i$. As in the proof of 
Lemma \ref{LemAccum}, we take $K_1$ so that $f^k(x_0)\in U_{i_1}$ 
for $k\geq K_1$ and set $k_n^{(1)} = min\{k\geq K_1; f^k(x_n)\notin U_{i_1}\}$. 
If $m_n < k_n^{(1)}$ for infinitely many $n$, 
then $p = \lim f^{m_n}(x_n)\in\ovl{U_{i_1}}$, a contradiction. 
Thus $m_n\geq k_n^{(1)}$ for large $n$. We may assume $f^{k_n^{(1)}}(x_n)$ 
tends to some $y^{(1)}\in W^u(\Lambda_{i_1})\setminus\Lambda$. 
Suppose $y^{(1)}\in W^s(\Lambda_{i_2})$. 
Since $y^{(1)}\in A(C_i)$, we have $i_2\in I_i$. 

Now take $K_2$ so that $f^k(y^{(1)})\in U_{i_2}$ for $k\geq K_2$ and 
$k_n^{(2)} = min\{k\geq k_n^{(1)}+K_2; f^k(x_n)\notin U_{i_2}\}$. 
If $k_n^{(1)}\leq m_n\leq k_n^{(1)}+K_2$ for infinitely many $n$, there exists 
$j\leq K_2$ so that $m_n = k_n^{(1)}+j$ for infinitely many $n$. Then we have 
$$
p = \lim f^{m_n}(x_n) = \lim f^{k_n^{(1)}+j}(x_n) = f^j(y^{(1)})\in W^u(\Lambda_{i_1}).
$$
Otherwise, we have $m_n\geq k_n^{(1)}+K_2$ for large $n$. 
We may assume $f^{k_n^{(2)}}(x_n)\to y^{(2)}\in W^u(\Lambda_{i_2})\setminus\Lambda$. 
If $m_n < k_n^{(2)}$ for infinitely many $n$, 
then $p = \lim f^{m_n}(x_n)\in\ovl{U_{i_2}}$, a contradiction. 
Thus $m_n\geq k_n^{(2)}$ for large $n$. 

Repeating this argument, we eventually meet $\Lambda_i$. That is, 
there exist $\ell$ and $i_j\in I_i, 1\leq j\leq\ell$ such that 
$$
\Lambda_{i_1}\succ\Lambda_{i_2}\succ\cdots\succ\Lambda_{i_{\ell}} = \Lambda_i. 
$$
Suppose $k_n^{(\ell)} < \infty$ for infinitely many $n$. Then, further repeating 
this argument, we must meet $\Lambda_i$ again. That is, there exists a sequence : 
$$
\Lambda_{i_1}\succ\cdots\succ\Lambda_{i_{\ell}} = \Lambda_i\succ\cdots\succ\Lambda_i. 
$$
This contradicts Theorem \ref{ThmNoCycle}. Thus we conclude that, for large $n$, 
$k_n^{(\ell)} = \infty$ and $f^k(x_n)\in U_{i_{\ell}}$ 
for $k\geq k_n^{(\ell-1)}+K_{\ell}$. Since $p\notin\ovl{U_{i_{\ell}}}$, 
we may conclude $k_n^{(\ell-1)}\leq m_n\leq k_n^{(\ell-1)}+ K_{\ell}$ for 
large $n$. Then, there exists $j\leq K_{\ell}$ such that $m_n = k_n^{(\ell-1)}+j$ 
for infinitely many $n$ and 
$$
p = \lim f^{m_n}(x_n) = \lim f^{k_n^{(\ell-1)}+j}(x_n) 
= f^j(y^{(\ell-1)})\in W^u(\Lambda_{i_{\ell-1}}). 
$$ 
Thus we conclude 
$A(C_i)\subset\cup_{j\in I_i} \left(W^u(\Lambda_j)\cap (J_p\times\CC)\right)$. 

Let us consider the case $i = 0$. Let $p\in A(C_0)$. The argument as above works 
as long as $y^{(\ell)}\in W^s(\Lambda)$. Suppose $y^{(\ell)}\notin W^s(\Lambda)$. 
Then it belongs to $W^s(\Lambda_0)$. 
Since $W^s(\Lambda_0)$ is open, for large $n$, $f^{k_n^{(\ell)}}(x_n)$ is contained 
in a neighborhood of $y^{(\ell)}$ in $W^s(\Lambda_0)$. 

If $m_n\leq k_n^{(\ell)}$ for infinitely many $n$, then $p\in\ovl{U_{i_{\ell}}}$, 
a contradiction. Thus $m_n\geq k_n^{(\ell)}$ for large $n$. 
Now suppose the sequence $\{m_n-k_n^{(\ell)}\}$ is unbounded. 
By taking a subsequence tending to $\infty$, the sequence 
$f^{m_n}(x_n) = f^{m_n-k_n^{(\ell)}}\circ f^{k_n^{(\ell)}}(x_n)$ 
tends to $\infty$, which is a contradiction. 
Hence $0 < m_n - k_n^{(\ell)} < K_{{\ell}+1}$ for some $K_{{\ell}+1} > 1$. 
Then there exists $j\geq 1$ such that $m_n = k_n^{(\ell)}+j$ 
for infinitely many $n$ and 
$p = \lim f^{m_n}(x_n) = f^j(y^{(\ell)})\in W^u(\Lambda_{i_{\ell}})$. 
This completes the proof. 
\qed

\vspace{2mm}
As a corollary, we extend the equivalence (3) in Theorem \ref{ThmAptEqA} to a 
union of saddle basic sets. Let $I'$ be a subset of $I$ and 
set $\Lambda' = \cup_{j\in I'}\Lambda_j,\: C' = \cup_{j\in I'}C_j.$ 

\begin{prop}\label{PropAC'Ws}
$A(C')\subset W^u(\Lambda')\cap W^s(\Lambda')\Longleftrightarrow 
C' \mbox{ is closed}. $
\end{prop}
\pr  $(\Leftarrow)$  Suppose $C'$ is closed. Then, for any $k\in I'$, 
it follows $I_k\subset I'$ since $\ovl{C_k}\subset \ovl{C'} = C'$. 
By Theorem \ref{ThmEstimate}, for any $k\in I'$, we have 
$$
A(C_k) \subset (\cup_{j\in I_k}W^u(\Lambda_j))\cap (\cup_{j\in I_k}W^s(\Lambda_j))
       \subset W^u(\Lambda')\cap W^s(\Lambda'). 
$$
Hence we have $A(C')\subset W^u(\Lambda')\cap W^s(\Lambda')$. 

\vspace{2mm}
\noindent
$(\Rightarrow)$  If $C'$ is not closed, there exists $k\in I'$ and $i\notin I'$ 
such that $\ovl{C_k}\cap C_i\neq\emptyset$. Then, by Proposition \ref{PropACiWs}, 
it follows $A(C_k)\cap W^s(\Lambda_i)\neq\emptyset$. Thus we conclude that 
$A(C')\not\subset W^u(\Lambda')\cap W^s(\Lambda')$. This completes the proof. 
\qed

\vspace{2mm}
We do not know if the equality holds in Proposition \ref{PropAC'Ws}. 
After Proposition \ref{PropAC'Ws}, a partition $I = \sqcup_i\tilde{I}_i$ of $I$ 
such that $\cup_{j\in\tilde{I}_i} C_j$ is closed for any $i$ will give equality. 
Define the equivalence relation among $I_i$ generated by $I_i\cap I_j\neq\emptyset$. 
That is, $I_i\sim I_j$ if there exist $i = i_1, i_2,\cdots,i_k = j$ in $I$ 
such that $I_{i_{\ell}}\cap I_{i_{\ell + 1}}\neq\emptyset$ for $1\leq\ell\leq k-1$. 
Let $I = \sqcup_i\tilde{I}_i$ be the partition of $I$ such that 
each $\tilde{I}_i$ is the union of $I_j$ in the same equivalence class. 
Then $I_j\subset\tilde{I}_i$ if $j\in \tilde{I}_i$ and 
it follows that the set $\cup_{j\in\tilde{I}_i}C_j$ is closed for each $i$. 
By Proposition \ref{PropAC'Ws}, for each $i$, we have 
$$
\cup_{j\in\tilde{I}_i} A(C_j)\subset (\cup_{j\in\tilde{I}_i}W^u(\Lambda_j))\cap 
(\cup_{j\in\tilde{I}_i}W^s(\Lambda_j)). \eqno{(7)}
$$
Since $A(C_{J_p}) = \cup_i\cup_{j\in\tilde{I}_i} A(C_j)$, the equality 
holds in (7). Thus we have 

\begin{prop}
$
\cup_{j\in\tilde{I}_i} A(C_j) = (\cup_{j\in\tilde{I}_i}W^u(\Lambda_j))\cap 
(\cup_{j\in\tilde{I}_i}W^s(\Lambda_j)). 
$
\end{prop}

\subsection{Proof of Theorem \ref{ThmStable}}
We use the stability of hyperbolic sets under 
perturbation established in \cite{j2,dh1,dh2}. 
Consider a holomorphic family $\{f_a;a\in\DD\}$ of Axiom A polynomial skew products 
containing $f = f_0$. The holomorphic motion $\Psi(a,z,w) = (\varphi_a(z),\psi_a(z,w))$ 
of $J_2 = J_2(f)$ gives homeomorphisms 
$$
\varphi_a : J_p\to J_{p_a},\quad \psi_a(z,\cdot) : J_z(f)\to J_{\varphi_a(z)}(f_a). 
$$
We have a holomorphic motion 
$\hat{h}_a : \hat\Lambda := \hat\Lambda(f)\to\hat\Lambda_a := \hat\Lambda(f_a)$ 
of $\hat\Lambda(f)$. As for $\Lambda := \Lambda(f)$, we only have a surjective 
continuous map $h_a : \hat\Lambda\to \Lambda_a := \Lambda(f_a)$, which is induced from 
$\hat{h}_a$ and depends holomorphically on $a$. Hence, for any $z\in J_p$, the map 
$a\mapsto\Lambda_{a,\varphi_a(z)} := \Lambda_a\cap (\{\varphi_a(z)\}\times\CC)$ 
is continuous. 
Note that the same holds also for each saddle basic set $\Lambda_i$ of $f$ 
and that the map $a\mapsto J_{\varphi_a(z)}(f_a)$ is continuous. 

Suppose a critical point $x = (z,w)\in C_{J_p}$ of $f$ lies in $C_i$. 
Then there exist $y = (u,v)\in\Lambda_i$ and $n\geq 0$ such that 
$x_n = f^n(x)\in W^s_{loc}(y)$. Then $U_{x_n}$, the connected component of the 
vertical slice of $(J_p\times\CC)\setminus J_2$ containing $x_n$, contains $y$ 
because $W^s_{loc}(y)$ lies in a vertical fiber. Set $z_a = \varphi_a(z)$ and 
let $x_a = (z_a,w_a)$ be a nearby critical point of $f_a$. 
If $a$ is small, $x_{a,n} = f_a^n(x_a)$ is close to $f^n(x)$. 
By the continuity of the maps $a\mapsto J_{z_a}(f_a)$ and $a\mapsto\Lambda_{a,z_a}$, 
the point $y_a = h_a(\hat{y})$ for any prehistory $\hat{y}$ of $y$ lies 
in $U_{x_{a,n}}\cap\Lambda_{a,i}$. Hence it follows $x_a\in W^s(\Lambda_{a,i})$ 
for small $a$. Theorem \ref{ThmAxiomA} says that the postcritical set is disjoint 
from the second Julia set. By the above continuity, we conclude 
that $x_a\in W^s(\Lambda_{a,i})$ holds as long as the holomorphic motion exists. 

Now suppose $A_{cc}(C_{J_p}) = A_{pt}(C_{J_p})$ for $f = f_0$. 
For any $C_a\in\C(C_{J_{p_a}})$, there exists a connected component $J_a$ 
of $J_{p_a}$ such that $C_a$ is a connected component of $C_{J_a}$. 
Then $J_a = \varphi_a(J)$ for some connected component $J$ of $J_p$ and 
there exists $C\in\C(C_{J_p})$ close to $C_a$, which is a connected component 
of $C_J$. By Theorem \ref{ThmAptEqAcc}, we have $C\subset C_i$ for some $i$. 
There exists $n$ such that, for any $x\in C$, $f^n(x)\in W^s_{loc}(y)$ 
for some $y\in\Lambda_i$. Thus we can apply the above argument 
simultaneously to any $x\in C$ and we conclude that 
$C_a\subset C_{a,i} := C_{J_{p_a}}\cap W^s(\Lambda_{a,i})$, 
hence $A_{cc}(C_{J_p}) = A_{pt}(C_{J_p})$ holds for $f = f_a$ as long as the 
holomorphic motion exists. 

Next suppose $A(C_{J_p}) = A_{pt}(C_{J_p})$ for $f = f_0$. 
By Theorem \ref{ThmAptEqA}, $A(C_i) = \Lambda_i$ holds for any $i$. 
There exists $n$ such that, for any $x\in C_i$, $f^n(x)\in W^s_{loc}(y)$ 
for some $y\in\Lambda_i$. Then nearby critical point of $f_a$ for small $a$ 
belongs to $W^s(\Lambda_{a,i})$. 
By Theorem \ref{ThmAptEqA}, all $C_i$ are closed, hence mutually 
disjoint compact sets. Thus so are $C_{a,i}$ and 
$A(C_{J_p}) = A_{pt}(C_{J_p})$ holds also for $f = f_a$ for small $a$. 
By Theorem \ref{ThmAxiomA}, this holds as long as the holomorphic motion exists. 

Since any maps in the same hyperbolic component are connected by a chain of disks 
where the holomorphic motion exists, we get the conclusion. \qed

\subsection{Proof of Theorem \ref{ThmAccEqA}}
First we estimate the component-wise accumulation set $A_{cc}(C_{J_p})$. 
The key tool is a result due to Qiu \& Yin \cite{qy} that any non-point 
component of a disconnected polynomial Julia set is preperiodic. 
Although a complete characterization of the set $A_{cc}(C_{J_p})$ is not 
known, our estimate is enough to prove Theorem \ref{ThmAccEqA}. 
Let $J_{per}$ be the union of the periodic non-point components of $J_p$. 
Then $J_{per} = J_p$ if $J_p$ is connected and $J_{per} = \emptyset$ 
if $J_p$ is a Cantor set. 

\begin{prop}\label{PropAcc}
$A_{cc}(C_{J_p})\subset\Lambda\cup (W^u(\Lambda)\cap (J_{per}\times\CC))$. 
\end{prop}
\pr  For any $C\in\C(C_{J_p})$, there exists a connected component $J_0$ 
of $J_p$ such that $C$ is a connected component of $C_{J_0}$. 
If $J_0$ is a point component, then $C_{J_0}$ consists of finitely many 
points, hence $A(C_{J_0})\subset\Lambda$. Otherwise, it is a 
preimage of a periodic component $J_1$ of $J_p$, i.e. 
$p^k(J_0) = J_1\subset J_{per}$ for some $k$. 
See \cite{qy}, Theorem in Section 5. 
Let $J_{-k}$ be the (finite) union of the connected components of 
$p^{-k}(J_{per})$. 
We have only to show that $\ovl{\cup_{k\geq 0}A(C_{J_{-k}})}\subset 
W^u(\Lambda)\cap (J_{per}\times\CC)$. By Theorem \ref{ThmAxiomA}, 
the set $X = \ovl{\cup_{k\geq 0} f^k(C_{J_{-k}})}$ is 
disjoint from $J_2$. Then Proposition 3.3 in \cite{dh1} says 
$A(X)\subset W^u(\Lambda)\cap (J_p\times\CC)$. 
Evidently we have $A(X)\subset J_{per}\times\CC$. Thus it follows 
$$
\ovl{\cup_{k\geq 0}A(C_{J_{-k}})} = \ovl{\cup_{k\geq 0}A(f^k(C_{J_{-k}}))} 
\subset A(X) \subset W^u(\Lambda)\cap (J_{per}\times\CC). 
$$
This completes the proof. 
\qed

\vspace{2mm}
We investigate the slice $W^u(\Lambda)_z$ of $W^u(\Lambda)$ at $z\in J_p$. 

\begin{lem}\label{LemUnstableNoEmpty}
Suppose $\Lambda\neq\emptyset$. Then $W^u(\Lambda)_z \neq\emptyset$ for 
any $z\in J_p$. 
\end{lem}
\pr  If $\Lambda\neq\emptyset$, there exists a saddle periodic point 
$x = (z_0,w_0)\in\Lambda$. Let $\hat{x}\in\hat{\Lambda}$ be any one of the 
prehistories of $x$. Recall that the set $\{p^{-k}(z);k\geq 0\}$ is dense 
in $J_p$ for any $z\in J_p$. Then, for any $z\in J_p$, there exist sequences 
$n_k\nearrow\infty$ and $z_{-k}\in p^{-n_k}(z)$ tending to $z_0$. 
Since $W^u_{loc}(\hat{x})$ is transversal to the vertical fiber, 
for large $k$, there exists $w_{-k}$ such that 
$x_{-k} := (z_{-k},w_{-k})\in W^u_{loc}(\hat{x})\cap (J_p\times\CC)$. 
Then $(z,w) := f^{n_k}(x_{-k})\in W^u(\Lambda)\cap (J_p\times\CC)$. 
This completes the proof. 
\qed

\vspace{2mm}
By virtue of Proposition \ref{PropAcc}, we get the following. 

\begin{prop}\label{PropLambda}
Suppose $J_p$ is disconnected and $\Lambda\neq\emptyset$. 
If $A_{cc}(C_{J_p}) = A(C_{J_p})$, then $W^u(\Lambda)_z = \Lambda_z$ 
for $z\in J_p\setminus J_{per}$ and 
$\Lambda_z\neq\emptyset$ for any $z\in J_p$. 
\end{prop}
\pr  Suppose $A_{cc}(C_{J_p}) = A(C_{J_p})$. Then, from Proposition 
\ref{PropAcc}, we have $W^u(\Lambda)\cap (J_p\times\CC)\subset 
\Lambda\cup (W^u(\Lambda)\cap (J_{per}\times\CC))$. Hence it follows 
$W^u(\Lambda)_z = \Lambda_z$ for $z\in J_p\setminus J_{per}$, 
which is not empty by Lemma \ref{LemUnstableNoEmpty}.

Note that $J_p\setminus J_{per}$ is dense in $J_{per}$. In fact, otherwise, 
a point in $J_{per}$ has a neighborhood $U$ disjoint from $J_p\setminus J_{per}$. 
Then $\cup_{n\geq 0}\, p^n(U)$ never intersects $J_p\setminus J_{per}$, 
a contradiction. Thus, for any $z\in J_{per}$, there exists 
$z_n\in J_p\setminus J_{per}$ tending to $z$. Take $(z_n,w_n)\in \Lambda$. 
Since $\{w_n\}$ is a bounded set, there is a subsequence $\{w_{n_k}\}$ converging 
to a point $w$. Then $(z,w) = \lim_{k\to\infty}(z_{n_k},w_{n_k})\in\Lambda$. 
This completes the proof. 
\qed

\vspace{2mm}
Let $d_H(\cdot,\cdot)$ be the Hausdorff distance of compact sets in $\CC$. 
Set $\DD(z_0,r) = \{z\in\CC;|z - z_0| < r\}$. 

\begin{prop}\label{PropConti}
Suppose $J_p$ is disconnected, $\Lambda\neq\emptyset$ and 
$W^u(\Lambda)_z = \Lambda_z$ for any $z\in J_p\setminus J_{per}$. Let $z_0$ be 
any point in $J_p$. Then, for any $\eps > 0$, there exists $\delta > 0$ such that 
$$
d_H(\Lambda_z,\Lambda_{z_0}) < \eps\quad \mbox{ for any }\quad
z\in\DD(z_0,\delta)\cap (J_p\setminus J_{per}).               \eqno{(8)} 
$$
In particular, the map $z\mapsto\Lambda_z$, restricted to $J_p\setminus J_{per}$, 
is continuous. 
\end{prop}
\pr  We prove the proposition by contradiction. Suppose there exist an $\eps > 0$ 
and a sequence $z_n\in J_p\setminus J_{per}$ tending to $z_0\in J_p$ 
with $d_H(\Lambda_{z_n},\Lambda_{z_0})\geq\eps$. Since $\Lambda$ is closed, 
the map $z\mapsto\Lambda_z$ is upper semicontinuous on $J_p$. Thus it follows
that $\Lambda_{z_0}$ does not lie in the $\eps$-neighborhood of $\Lambda_{z_n}$, 
that is, there exists $w_n\in\Lambda_{z_0}$ such that 
$\Lambda_{z_n}\cap\DD(w_n,\eps) = \emptyset$. 
We may assume $w_n\to w_0\in\Lambda_{z_0}$. 
Then $\Lambda_{z_n}\cap\DD(w_0,\eps/2) = \emptyset$ for large $n$. 
Since $W^u_{loc}(\hat{x})$ is transversal to the fiber for any prehistory 
$\hat{x}$ of $x = (z_0,w_0)\in\Lambda$, we have 
$W^u(\Lambda)_{z_n}\cap\DD(w_0,\eps/2)\neq\emptyset$ for large $n$. 
This is a contradiction since $W^u(\Lambda)_{z_n} = \Lambda_{z_n}$. 
\qed

\vspace{2mm}
Now Theorem \ref{ThmAccEqA} follows from the next theorem. 

\begin{thm}\label{ThmConti}
Suppose $J_p$ is disconnected. If $W^u(\Lambda)_z = \Lambda_z$ for 
any $z\in J_p\setminus J_{per}$, then the map $z\mapsto\Lambda_z$ is 
continuous in $J_p$. In particular, if $A_{cc}(C_{J_p}) = A(C_{J_p})$, 
then $A_{pt}(C_{J_p}) = A_{cc}(C_{J_p}) = A(C_{J_p})$. 
\end{thm}
\pr  Theorem is true if $\Lambda = \emptyset$. Hence we may assume 
$\Lambda\neq\emptyset$. Take a point $z_0$ in $J_p$. We have only to 
show (8) also for $z\in\DD(z_0,\delta)\cap J_{per}$. 
For any $\eps > 0$, take $\delta > 0$ as in Proposition \ref{PropConti}. 
Now, for any $z\in\DD(z_0,\delta)\cap J_{per}$, 
there exists a sequence $\{z_n\}$ in $J_p\setminus J_{per}$ tending to $z$. 
Applying Proposition \ref{PropConti} to $z_0 = z$, it follows 
$\Lambda_{z_n}\to\Lambda_z$. Note that $z_n\in\DD(z_0,\delta)$ for large $n$. 
Again by Proposition \ref{PropConti}, 
we have $d_H(\Lambda_{z_n},\Lambda_{z_0}) < \eps$ for such $n$, 
hence $d_H(\Lambda_z,\Lambda_{z_0})\leq\eps$. 
This implies the desired continuity on $J_{per}$. 
The last statement follows from Theorem C. This completes the proof. 
\qed

\subsection{Proof of Theorem \ref{ThmApt=AccNotA}}

The proof of Theorem \ref{ThmApt=AccNotA} is a higher degree analogue of that 
of Theorem 6.1 in \cite{dh1}. Consider the degree four polynomial skew product 
$f(z,w) = (p(z),q(z,w))$, where $q(z,w) = w^4 + 4(2-z)$ and $p$ is a 
{\it real biquadratic polynomial} of the form : 
$$
p(z) = p_{a,b}(z) = (z^2+a)^2+b,\quad (a,b)\in\RR^2.
$$

First we investigate the dynamics of $p$ so that we find maps with one 
critical point escaping while the others tending to an attracting cycle. 
Next we show that the map $f$ has the desired property. 
For our maps, the saddle set $\Lambda$ consists of a single saddle fixed point. 
Then exactly one critical point tends to this point while all 
the others escape to infinity. The hard part is to show that 
$D_{J_p}\cap J_2 = \emptyset$. Most of this subsection, 
Lemmas \ref{LemNumber} - \ref{LemVertExpand}, is devoted to control 
the behavior of the critical orbits. 

Here we give a brief summary on the external rays of a monic polynomial $p$ of degree $d$. 
See Milnor \cite{mil} for the details. Let $\varphi_p$ be the {\it B\"{o}ttcher coordinate} 
of $p$. It is a conformal map defined in a neighborhood of $\infty$ conjugating 
$p$ to the map $z\mapsto z^d$. 
Then the {\it external ray} $R_p(\theta)$ of $p$ with angle $\theta$ 
is defined by $R_p(\theta) = \varphi_p^{-1}(\{re^{2\pi i\theta};r > r_p\})$ for some 
$r_p\geq 1$ depending on $p$. It maps $R_p(\theta)$ to $R_p(d\theta)$. 
Hence, it is continued until it meets a critical point. 
If it is continued to $r > 1$ and the limit 
$z = \lim_{r\to 1}\varphi_p^{-1}(re^{2\pi i\theta})$ 
exists, it is said to land at $z$. A fundamental fact is as follows : 
for $\theta\in\QQ$, the $\theta$-ray $R_p(\theta)$ lands at some point 
$z\in J_p$ unless it meets a point in the backward orbit of a critical point. 
If the 0-ray lands, its landing point is a fixed point, which we call the 
$\beta$-fixed point of $p$ and denote by $\beta_p$. 

The map $p = p_{a,b}$ has three critical points : $z = 0,\:\pm\sqrt{-a}$ and 
two critical values : $p(0) = a^2+b,\: p(\pm\sqrt{-a}) = b$. Note that critical 
values are always real. The {\it connectedness locus} $\C$ of this family is 
the set of parameters $(a,b)$ in $\RR^2$ so that the Julia set 
$J(p_{a,b})$ is connected, or equivalently, all the critical points of $p_{a,b}$ 
have bounded orbits. It is described as follows. See Figure \ref{FigM}. 
The dark region indicates $\C$. 

\begin{lem}\label{LemConnLocus}
The boundary of $\C$ consists of the following three curves :

\noindent
$\displaystyle Per_1^+(1): p'(\beta_p) = 1 : (a,b) = (\frac{1}{2t} - \frac{t^2}{4},
\frac{t}{2} - \frac{1}{4t^2}),\: \frac{1}{\sqrt[3]{4}}\leq t\leq\sqrt[3]{4}$, 

\vspace{1mm}
\noindent
$\displaystyle Preper_{(1)1}: p(0) = \beta_p : b = -a^2 + \sqrt{-2a},\: 
- 2 \leq a\leq - \frac{\sqrt[3]{2}}{4}$, \\
$\displaystyle Preper_{(2)1}: p(b) = \beta_p = - b : a = -b^2 + \sqrt{-2b},\: 
- 2 \leq b\leq - \frac{\sqrt[3]{2}}{4}$. 
\end{lem}
\pr  Since the critical values are real, we consider the dynamics of $p$ only 
on the real axis. For large $b$, the graph of $y = p(x)$ sits above the diagonal 
line $y = x$, hence the orbits of all real points tend to $+\infty$. 
Decreasing $b$, we meet a parameter at which the graph of $y = p(x)$ is tangent 
to the line $y = x$ at $\beta_p$. This parameter lies on $Per_1^+(1)$, where 
$\beta_p$ is a parabolic fixed point with multiplier one. 
Below this locus, the map $p$ has at least two real fixed points. The point 
$\beta_p$ is the largest one. Then $p\in\C$ if and only if 
$K_p\cap\RR = [-\beta_p,\beta_p]$. 
If $a < 0$, then $p$ has the local maximum $p(0)$ and the local minimum $b$. 
Hence $p\in\C$ if and only if $-\beta_p\leq b$ and $p(0)\leq\beta_p$. 
In case $a\geq 0$, since $p$ has a unique local minimum $p(0) = a^2 + b\geq b$, 
$p\in\C$ if and only if $b\geq -\beta_P$. 

The locus $Preper_{(1)1} : p(0) = \beta_p$ is included in the set : 
$\{(a^2 + b)^2 + a\}^2 + b = a^2 + b,$ that is, $(a^2+b)^2 + a = \pm a$. 
Since $(a^2+b)^2 + a = a$ implies $p(0) = 0$, $Preper_{(1)1}$ is included 
in $(a^2+b)^2 = -2a$ i.e. $b = -a^2 \pm\sqrt{-2a}$. 
Since $p(0) = a^2 + b = \beta_p > 0$, we conclude that $Preper_{(1)1}$ is 
written as $b = -a^2 + \sqrt{-2a}$. 

The locus $Preper_{(2)1} : b = - \beta_p$ or $p(b) = - b$ is included in 
$(b^2+a)^2 + 2b = 0$, that is, $b\leq 0$ and $a = -b^2 \pm\sqrt{-2b}$. 
On the curve $a = -b^2 - \sqrt{-2b}$, we have 
\begin{eqnarray*}
p(x) - x &=& x^4 - 2(b^2+\sqrt{-2b})x^2 - x + b^4 + 2b^2\sqrt{-2b} - b\\
         &=& (x + b)(x^3 - bx^2 - (b^2+2\sqrt{-2b})x + b^3 + 2b\sqrt{-2b} - 1)\\
         &=:& (x+b)g(x),
\end{eqnarray*}
and $g(-b) = 4b\sqrt{-2b} - 1 < 0$. Thus $p$ has a fixed point larger than $-b$. 
That is, $\beta_p > - b$, hence $Preper_{(2)1}$ is written by 
$a = -b^2 + \sqrt{-2b}$. This completes the proof. 
\qed 

\begin{figure}[h]
\noindent
\begin{minipage}{70mm}
\centering{
\unitlength 0.1in
%
}
\caption{Parameter space}
\label{FigM}
\end{minipage}
\vskip -71mm\hskip 70mm
\begin{minipage}{70mm}
\centering{
\unitlength 0.1in
%
%
}
\caption{Graph of $p_{-2,-2}$}
\label{FigR1}
\end{minipage}
\end{figure}

We consider the subfamily $Preper_{(1)1}$. 
Note that $p_{-2,-2}$ is the second iterate of the Chebyshev polynomial $z^2-2$ 
of degree two. 

\begin{lem}
There exists a sequence $(a_n,b_n)\in Preper_{(1)1}$ tending to $(-2,-2)$ such that 
$\sqrt{-a_n}$ is a superattracting periodic point of $\tilde{p}_n := p_{a_n,b_n}$ 
of period $n$. 
\end{lem}
\pr  For $(a,b)\in Preper_{(1)1}$, take an increasing sequence 
$x_n\in p^{-n}(\sqrt{-a})$ in $\RR$ tending to $\beta_p$. Then there exists a 
unique $n$ such that $x_{n-1}\leq p^2(\sqrt{-a}) < x_n$. Moving along $Preper_{(1)1}$, 
we get a parameter $(a_n,b_n)$ satisfying $p^2(\sqrt{-a}) = x_{n-1}$, 
that is $p^{n+1}(\sqrt{-a}) = \sqrt{-a}$. 
If the parameter approaches $(-2,-2)$, $p^2(\sqrt{-a})$ tends to $\beta_p$, 
hence it becomes larger than $x_n$ for arbitrarily large $n$. 
This completes the proof. 
\qed

By Lemma \ref{LemConnLocus}, a small perturbation from $Preper_{(1)1}$ 
leaves the connectedness locus. Thus, we get a small perturbation $p_n$ 
of $\tilde{p}_n$ outside the connectedness locus. 
For $p_n$, $0$ is escaping and the other critical points $\pm\sqrt{-a}$ belong 
to the basin of a superattracting $n$-cycle of $p_n$. 
Thus $p = p_n$ satisfies the property (a) in Theorem \ref{ThmApt=AccNotA}. 
See Figures \ref{FigR1} - \ref{FigPerturb}. 

\begin{figure}[h]
\begin{minipage}{70mm}
\centering{
\unitlength 0.1in
%
}
\caption{Graph of $\tilde{p}_n$}
\label{FigNonPerturb}
\end{minipage}
\vskip -71mm\hskip 70mm
\begin{minipage}{70mm}
\centering{
\unitlength 0.1in
%
%
}
\caption{Graph of $p_n$}
\label{FigPerturb}
\end{minipage}
\end{figure}

By upper-semicontinuity of the filled-in Julia set, we have the following. 

\begin{lem}\label{LemKset}
There exists a sequence $\eps_n$ tending to $0$ such that \\
$K_{p_n}\subset [-5/2,5/2]\times [-\eps_n,\eps_n]$. 
\end{lem}

Now we will control the behaviour of the critical orbits. 

\begin{lem}\label{LemFiberKset}
There exists $n_0$ such that for any $n\geq n_0$ and for any $z\in K_{p_n}$, we have \\
$(a)$ $q_z(\CC\setminus\DD(0,5/2))\subset\CC\setminus\DD(0,35/2)$,\\
$(b)$ $K_z \subset \DD(0,5/2)$. 
\end{lem}
\pr  By Lemma \ref{LemKset}, there exists $n_0$ such that $\eps_n\leq 1/4$ 
for $n\geq n_0$. Then $|2-z|\leq 5$ and, if $|w|\geq 5/2$, we have 
$$
|q_z(w)| \geq |w|^4 - 4|2-z| \geq |w|(|w|^3 - \frac{20}{|w|})
         \geq (\frac{125}{8} - 8)|w| \geq 7|w| \geq \frac{35}{2}.
$$
This completes the proof. 
\qed

\vspace{2mm}
Set $f_n(z,w) = (p_n(z),q(z,w))$ and $f_n^k(z,w) = (p_n^k(z),Q_{n,z}^k(w))$. 

\begin{lem}\label{LemVertExpA}
For $n\geq n_0$, $D_{A_{p_n}}\cap J_{A_{p_n}} = \emptyset$. 
\end{lem}
\pr  Note that $A_{p_n}\subset\RR$ and $C_{A_{p_n}} = A_{p_n}\times\{0\}$. 
Take $x\in A_{p_n}$ and set $(x_k,y_k) = f_n^k(x,0)$. Then $x_k = p_n^k(x)$. 
Recall that an attracting cycle contains the critical value 
$b = p_n(\pm\sqrt{-a})$ in the cycle of its immediate basins. 
Since $b = b_n < 0$, the basin $U$ containing $b$ satisfies 
$U\cap\RR\subset\{\Re\, z < 0\}$. Hence $x_j < 0$ for some $0\leq j\leq n-1$. 
Then $y_{j+1} = y_j^4 + 4(2-x_j)\geq 8$. By Lemma \ref{LemFiberKset}, 
$y_{j+1}\notin K_{x_{j+1}}$, hence $y_k\notin K_{x_k}$ for any $k$. 
Thus we conclude $D_{A_{p_n}}\cap J_{A_{p_n}} = \emptyset$. 
This completes the proof. 
\qed

\vspace{2mm}
In the following lemma, the constant $r$ is fixed. Later we assume $r < 7/128$. 
Set $B(z_0,r) = \{z\in\CC;|z - z_0|\leq r\}$. 

\begin{lem}\label{LemNumber}
For any fixed $r > 0$, there exists $N > 0$ and $n_1\geq n_0$ such that for 
any $n\geq n_1$ and for any $z\in J_{p_n}\setminus B(2,r)$, there exists 
$0\leq j < N$ with $\Re\, z_j\leq 1$. 
\end{lem}
\pr  For $p_{\infty} := p_{-2,-2}$, the $\pm 1/4$-rays land at $0$ 
and the $\pm 3/16$-rays land at $x_{\infty} = \sqrt{2-\sqrt{2}} = 0.7655\ldots$, 
a preimage of $0$. Thus the set of landing points $z\in J_{p_{\infty}}$ 
of external rays with angles in $[3/16,13/16]$ sits in $\Re\, z \leq 1$. 
We will show the same property also for $p_n$ for large $n$. 

First we consider maps $\tilde{p}$ on $Preper_{(1)1}$. 
The 0-ray lands at $\beta_{\tilde{p}}$. 
By the symmetry of the filled-in Julia set, the $1/4$-ray and the $-1/4$-ray 
are respectively the positive and negative imaginary axes. 
As their preimages, the $\pm 3/16$-rays and their landing points depend 
continuously on $\tilde{p}\in Preper_{(1)1}$. By the upper-semicontinuity 
of the filled-in Julia sets, we conclude that the same property holds 
for $\tilde{p}_n$ for large $n$. 

Recall that $p_n$ is obtained by a small perturbation of $\tilde{p}_n$ 
outside the connectedness locus. Let $p$ be such a perturbation 
of $\tilde{p}\in Preper_{(1)1}$. It is easy to see that, for $p$, 
the 0-ray also lands at $\beta_p$. 
Again by the symmetry, the positive and negative imaginary axes respectively 
form the $1/4$- and $-1/4$-rays for $p$ and meet at the escaping critical 
point $0$. As their preimages, the $3/16$- and $-3/16$-rays in the upper and 
lower half plane respectively meet at a preimage of $0$ and depend continuously 
on $p$. By the same way as above, we conclude that the set of points 
$z\in J_{p_n}$ with external angles in $[3/16,13/16]$ sits in $\Re\, z \leq 1$ 
for large $n$. 

Note that the part of $K_{p_{\infty}}$ on the right side of the $\pm 3/4^j$-rays 
sits in $B(2,r)$ for large $j$. 
The same holds for $p$ sufficiently close to $p_{\infty}$, 
since the $\pm 3/4^j$-rays depend continuously on the parameter. 

Now define a map $m :\RR/\ZZ\to\RR/\ZZ$ by $m(t) = 4t$. 
Then $m^{-1}((3/16,13/16))$ contains $(3/64,13/64)\cup (51/64,61/64)$ and 
the sum of the preimages \\ 
$\cup_{j\geq 2}\{(3/4^j,13/4^j)\cup(-13/4^j,-3/4^j)\}$ covers $\RR/\ZZ$ 
except the angle $\{0\}$. The part of $K_{p_n}$ on the right side 
of the $\pm 3/4^j$-rays is included in $B(2,r)$ for large $j$. 
This completes the proof. 
\qed

\vspace{2mm}
Note that, the $\beta$-fixed point $\beta_n := \beta_{p_n}$ forms a single 
point component of $J_{p_n}$ because the $\pm 3/4^j$-rays for $j\geq 2$ separate 
$\beta_n$ from any other points in $J_p$. 

\begin{lem}\label{LemEscape}
Let $N$ be given in Lemma \ref{LemNumber}. 
There exists $n_2\geq n_1$ and $\delta > 0$ such that 
$$
Q_{n,z}^N(\{|\Im\, w| < \delta\})\cap\DD(0,5/2) = \emptyset
$$
for any $z\in J_{p_n}\setminus B(2,r)$ and for any $n\geq n_2$. 
\end{lem}
\pr  Take $z\in J_{p_n}\setminus B(2,r)$ and set $(z_k,w_k) = f_n^k(z,w)$. 
Let $j \leq N-1$ be the minimum of $k$ with $\Re\, z_k\leq 1$, which is assured by 
Lemma \ref{LemNumber}. Suppose $|\Re\, w_k| > 5/2$ for some $k\leq j$. Then 
$w_k\notin \DD(0,5/2)$, hence by Lemma \ref{LemFiberKset}, $w_N\notin\DD(0,5/2)$. 
Thus we may assume $|\Re\, w_k|\leq 5/2$ for $k\leq j$. 
We take $\delta$ and $n_2$ so that 
$$
64^N(\delta + \frac{4\eps_n}{63}) < \frac{\sqrt{6}}{10}, 
$$
for any $n\geq n_2$. 

Set $w = u + iv$ and $w_k = u_k + iv_k$. We show, by induction on $k$, 
$\displaystyle |v_k| < \frac{\sqrt{6}}{10}$ for $k\leq j$. 
Suppose it is true for $k = m$. Then it follows 
\begin{eqnarray*}
|v_{m+1}| &=& |4(u_m^2-v_m^2)u_mv_m - 4\Im\, z_m|\\
          &\leq& 4(u_m^2+v_m^2)|u_mv_m| + 4\eps_n\\
          &\leq& 4((\frac{5}{2})^2 + v_m^2)\frac{5}{2}|v_m| + 4\eps_n. 
\end{eqnarray*}
By the induction hypothesis, we have 
$\displaystyle 4((\frac{5}{2})^2 + v_m^2)\frac{5}{2} < \frac{631}{10} < 64$. 
Thus we have 
$|v_{m+1}| < 64|v_m| + 4\eps_n$ and it follows 
$$
|v_{m+1}| < 64^{m+1}|v| + 4\eps_n\sum_{j=0}^m64^j 
          < 64^N(\delta + \frac{4\eps_n}{63}) 
           <  \frac{\sqrt{6}}{10}. 
$$
Thus the case $k = m+1$ holds. 

Since $\Re\, z_j \leq 1$ from the choice of $j$, we have 
\begin{eqnarray*}
u_{j+1} &=& (u_j^2-v_j^2)^2 - 4u_j^2v_j^2 + 4(2-\Re\, z_j)\\
         &\geq& 4(2-\Re\, z_j) - 4u_j^2v_j^2\\
         &>& 4 - 4\cdot\frac{25}{4}\cdot\frac{6}{100} = \frac{5}{2}. 
\end{eqnarray*}
Thus $w_{j+1}\notin\DD(0,5/2)$, hence $w_N\notin\DD(0,5/2)$ 
by Lemma \ref{LemFiberKset}, which proves the lemma. 
\qed

\vspace{2mm}
Set $S(0,1/4,r) = \{x+yi\in\CC; |x|\leq 1/4, |y|\leq r\}$. 

\begin{lem}\label{LemContract}
Suppose $r < 7/128$. Then, for any $\delta' < 1/4$, there exists $n_3\geq n_2$ 
so that for all $n\geq n_3$ and $z\in J_{p_n}\cap B(2,r)$, 
we have $q_z(S(0,1/4,\delta'))\subset int\, S(0,1/4,\delta')$. 
\end{lem}
\pr  Fix $\delta' < 1/4$. Take $n_3$ so that $\eps_n\leq \delta'/8$ for $n\geq n_3$. 
Take $w\in S(0,1/4,\delta')$ and $z\in J_{p_n}\cap B(2,r)$. Then 
\begin{eqnarray*}
|v_1| &\leq& 4|u^2-v^2||uv| + 4\eps_n\\
      &\leq& 4(\frac{1}{16}+\delta'^2)\frac{1}{4}\delta' + 4\eps_n < \delta'. \\
|u_1| &\leq& (u^2-v^2)^2 + 4u^2v^2 + 4(2 - \Re\, z)\\
      &\leq& (\frac{1}{4^2}+\delta'^2)^2 + \frac{4\delta'^2}{4^2} + 4r\\
      &\leq& \frac{1}{8^2} + \frac{1}{4^3} + 4r < \frac{1}{4}.
\end{eqnarray*}
This completes the proof. 
\qed

\vspace{2mm}
Let $\delta$ and $n_2$ be given in Lemma \ref{LemEscape} and let $n_3$ be given 
in Lemma \ref{LemContract}. 

\begin{lem}\label{LemCritDisjoint}
For $n\geq n_3$, we have \\
$(a)$ $S(0,1/4,\delta)\cap K_z = \emptyset$ 
for any $z\in J_{p_n}\setminus\{\beta_n\}$,\\
$(b)$ $S(0,1/4,\delta)\cap J_{\beta_n} = \emptyset$, \\
$(c)$ $J_{p_n}\times S(0,1/4,\delta) \subset (J_{p_n}\times\CC)\setminus J_2$. 
\end{lem}
\pr  $(a)$\:  Any $z\in J_{p_n}\setminus\{\beta_n\}$ leaves $B(2,r)$ under finite 
iterations of $p_n$ because $\beta_n$ is a repelling fixed point of $p_n$. 
Set $m = min \{k;z_k := p_n^k(z)\notin B(2,r)\}$. 
By Lemma \ref{LemContract}, it follows 
$$
Q_{n,z}^m(S(0,1/4,\delta))\subset int\, S(0,1/4,\delta)\subset\{|\Im\, w| < \delta\}.
$$
Then applying Lemma \ref{LemEscape} to $z = z_m$, we have 
$$
Q_{n,z}^{m+N}(S(0,1/4,\delta))\cap K_{z_{m+N}} = \emptyset.
$$
By the invariance of $K$, we have $S(0,1/4,\delta)\cap K_z = \emptyset$. 

\noindent
$(b)$\:  Since $p_n(\beta_n) = \beta_n$, Lemma \ref{LemContract} implies 
$q_{\beta_n}(S(0,1/4,\delta))\subset int\, S(0,1/4,\delta)$. 
Then $S(0,1/4,\delta)\subset int\, K_{\beta_n}$ and 
we conclude $S(0,1/4,\delta)\cap J_{\beta_n} = \emptyset$. 

\noindent
$(c)$\:  It follows from $(a), (b)$. 
\qed

\begin{lem}\label{LemVertExpand}
For $n\geq n_3$, $D_{J_{p_n}}\cap J_2 = \emptyset$. 
\end{lem}
\pr  Set $(z_k,w_k) = f_n^k(z,0)$ for $z\in J_{p_n}$. 
If $z = \beta_n$, Lemma \ref{LemContract} says that $w_k\in S(0,1/4,\delta)$ 
for any $k\geq 0$ and by Lemma \ref{LemCritDisjoint} $(b)$, $\{(z_k,w_k)\}$ 
is uniformly bounded away from $J_2$. 

For $z \in J_{p_n}\setminus\{\beta_n\}$, set $m$ as above. 
If $m > 0$, Lemma \ref{LemContract} says $w_k\in int\, S(0,1/4,\delta)$ 
for $k\leq m$ and $\{(z_k,w_k);k\leq m\}$ is 
uniformly bounded away from $J_2$. Applying Lemma \ref{LemEscape} 
to $z = z_m$,  $\{(z_k,w_k);k\geq m + N\}$ is uniformly bounded away 
from $J_2$. The fact that $\{(z_k,w_k);m < k < m + N\}$ is uniformly 
bounded away from $J_2$ follows from the compactness of $J_{p_n}$ and 
the uniformity of $N$. If $m = 0$, i.e., $z\notin B(2,r)$, applying 
Lemma \ref{LemEscape} to $z$, we get the same conclusion. 
Thus we have shown that $D_{J_{p_n}}\cap J_2 = \emptyset$. 
\qed

\vspace{2mm}
Let $\delta$ and $n_2$ be given in Lemma \ref{LemEscape} and let $n_3$ be given 
in Lemma \ref{LemContract}. We may assume $\delta < 1/4$ and $n_3\geq n_2$. 

\begin{thm}\label{Thm1}
Set $f_n(z,w) = (p_n(z),q(z,w))$ as above. For $n\geq n_3$, $f_n$ 
is Axiom A and satisfies \\
$(a)$ $C_1 = \{(\beta_n,0)\},\quad C_0 = C_{J_{p_n}}\setminus C_1$. \\
$(b)$ $\Lambda = \Lambda_1 = \{(\beta_n,\alpha_n)\}$ is itself a basic set, 
where $\alpha_n$ is the attracting fixed point of $q_{\beta_n}$.\\
$(c)$ $J_z$ is disconnected for all $z\in J_{p_n}\setminus\{\beta_n\}$ 
and is a quasicircle for $z=\beta_n$. \\
$(d)$ $A_{pt}(C_{J_{p_n}}) = A_{cc}(C_{J_{p_n}})\neq A(C_{J_{p_n}})$. 
\end{thm}
\pr  By Theorem \ref{ThmAxiomA} and Lemmas \ref{LemVertExpA} 
and \ref{LemVertExpand}, $f_n$ for $n\geq n_3$ are Axiom A. 

\noindent
$(a)$\:  Lemma \ref{LemCritDisjoint}, $(a)$ says that 
$K\cap (C_{J_{p_n}}\setminus\{(\beta_n,0)\}) = \emptyset$. 
By Lemma \ref{LemContract}, it follows that $(\beta_n,0)\in K$. These imply $(a)$. 

\noindent
$(b)$\:  By Lemma \ref{LemContract}, $q_{\beta_n}$ has an attracting fixed point. 
From $(a)$, there is no saddle periodic point for $z\neq\beta_n$. 

\noindent
$(c)$\:  The case $z\in J_{p_n}\setminus\{\beta_n\}$ follows from $(a)$ and 
Proposition 2.3 in \cite{j3}. If $z=\beta_n$, $q_z(w) = w^4 + 4(2-\beta_n)$ has 
an attracting fixed point $\alpha_n$. Thus $J_{\beta_n}$ is a quasicircle 
as a quasiconformal image of the Julia set of $w\mapsto w^4$. 

\noindent
$(d)$\:  Since $C_{J_{p_n}} = J_{p_n}\times\{0\}$, any $C \in\C(C_{J_{p_n}})$ is 
of the form $J\times\{0\}$, where $J$ is a connected component of $J_{p_n}$. 
Recall that the component of $J_{p_n}$ containing $\beta_n$ is just 
a single point $\{\beta_n\}$ itself. Thus $(a)$ implies $C\subset C_0$ 
or $C\subset C_1$. Now by Theorem \ref{ThmAptEqAcc}, 
we conclude $A_{cc}(C_{J_{p_n}}) = A_{pt}(C_{J_{p_n}})$. 

Since $C_0$ is not closed in $C_{J_p}$, it follows from Theorem \ref{ThmAptEqA}, 
that $A(C_{J_{p_n}})\neq A_{pt}(C_{J_{p_n}})$, 
Thus we conclude $A_{pt}(C_{J_{p_n}}) = A_{cc}(C_{J_{p_n}})\neq A(C_{J_{p_n}})$. 
\qed

\end{document}